\documentclass[11pt]{article} 
\usepackage[utf8]{inputenc}
\usepackage[T1]{fontenc}
\usepackage{lmodern}
\usepackage[normalem]{ulem}
\usepackage[english]{babel}

\usepackage{amsmath,amsfonts,amssymb,amsthm}
\usepackage{marvosym}

\usepackage[mono=false]{libertine}
\usepackage[cmintegrals,libertine]{newtxmath}
\usepackage[cal=euler, scr=boondoxo]{mathalfa}

\useosf
\usepackage[a4paper,vmargin={3.5cm,3.5cm},hmargin={2.5cm,2.5cm}]{geometry}
\usepackage[font={small,sf}, labelfont={sf,bf}, margin=1cm]{caption}
\usepackage[pdftex,colorlinks=true]{hyperref}
\usepackage{color,graphicx}

\usepackage{stackrel}

\theoremstyle{plain}
\newtheorem{theorem}{Theorem}

\newtheorem{proposition}[theorem]{Proposition}
\newtheorem{lemma}[theorem]{Lemma}
\theoremstyle{definition}
\newtheorem{definition}[theorem]{Definition}
\newtheorem*{question}{Question}
\newtheorem{remark}[theorem]{Remark}

\newcommand{\expo}{  \mathsf{e}_{ \qseq}}
\newcommand{\gulp}{  \mathsf{g}_{ \qseq}}
\newcommand{\sexpo}{  \hat{\mathsf{e}}_{ \qseq}}
\newcommand{\sgulp}{  \hat{\mathsf{g}}_{ \qseq}}

\newcommand{\Map}{\mathfrak{M}}
\newcommand{\map}{\mathfrak{m}}

\newcommand{\rootface}{f_{ \mathrm{r}}}

\newcommand{\qseq}{\mathbf{q}}

\newcommand{\hatsubset}{\mathbin{\hat{\subset}}}
\newcommand{\hatSubset}{\mathbin{\hat{\Subset}}}

\begin{document}

\title{\bf Simple peeling of planar maps\\ with application to site percolation}
\author{\textsc{Timothy Budd}\footnote{Radboud University, Nijmegen, The Netherlands. Email: \href{mailto:t.budd@science.ru.nl}{t.budd@science.ru.nl}} \, and \textsc{Nicolas Curien}\footnote{Universit\'e Paris-Saclay and Institut Universitaire de France. E-mail: \href{mailto:nicolas.curien@gmail.com}{nicolas.curien@gmail.com}.}}
\date{\today}
\maketitle


\begin{abstract}
The peeling process, which describes a step-by-step exploration of a planar map, has been instrumental in addressing percolation problems on random infinite planar maps.
Bond and face percolation on maps with faces of arbitrary degree are conveniently studied via so-called lazy-peeling explorations.
During such explorations distinct vertices on the exploration contour may at later stage be identified, making the process less suited to the study of site percolation.
To tackle this situation and to explicitly identify site-percolation thresholds, we come back to the alternative ``simple'' peeling exploration of Angel and uncover deep relations with the lazy-peeling process. Along the way we define and study the random Boltzmann map of the half-plane with a simple boundary for an arbitrary critical weight sequence. 
Its construction is nontrivial especially in the ``dense regime'' where the half-planar random Boltzmann map does not possess an infinite simple core.
\end{abstract}

\vspace{1cm}

\begin{figure}[h]
\begin{center}
		\includegraphics[width=0.5\linewidth]{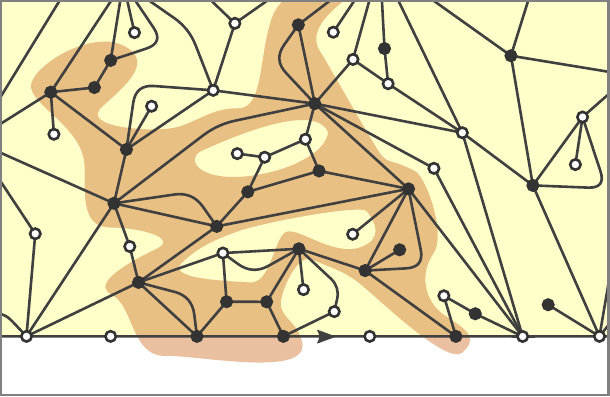}
		\caption{ \label{fig:sitepercolation} Illustration of Bernoulli site percolation on an infinite bipartite planar map of the half plane with simple boundary. The shaded region indicates the black cluster containing the origin of the root edge. }
	\end{center}
\end{figure}

\section{Introduction} 
Conceived by Watabiki \cite{Wat95} in the early 90's and formalized by Angel \cite{Ang03}, the \emph{peeling process} is a step-by-step exploration of random planar maps and a key tool for studying those random lattices. The peeling considered by Angel, based on maps with \emph{simple} boundaries is what we call in this work the \emph{simple peeling process}. It was analyzed in depth in the case of triangulations and quadrangulations (see e.g.~\cite{Ang03,Ang05,ACpercopeel,AR13,Ray13,MN13,CCsaw,BCsubdiffusive,CurKPZ,gwynne2017convergence,GM16a,bernardi2018percolation,BCKscalingpeeling}) but not for maps with high face-degrees due to the complicated transitions of the process and the lack of explicit enumerative formulas for maps with simple boundaries. An alternative peeling process closer to the initial construction of Watabiki, called here the \emph{lazy-peeling process}, was introduced by the first author in \cite{Bud15}. The lazy-peeling process has only two different ``topological types'' of transitions which have a universal form for all Boltzmann maps, making its analysis more transparent compared to the simple peeling. See \cite{CurStFlour} for a comprehensive summary. In this work, we come back to the simple peeling process on Boltzmann maps and uncover deep relations with the lazy-peeling process. As an application, we are able to give a closed formula for the site percolation threshold which seems out of reach using the lazy-peeling process only. On the way, we study Boltzmann maps with a simple boundary and their local limits, introducing in particular the half-planar Boltzmann maps with a simple boundary whose existence is new in the so-called dense regime.

\paragraph{Prerequisites on general Boltzmann maps.} In this work we focus on bipartite rooted planar maps in order to stick to the same framework as the reference \cite{CurStFlour}, from which we borrow notation and to which we often refer for background. Although computations will be more tedious, we expect our results and proofs to extend straightforwardly to the case of non-bipartite maps (see \cite{Bud15} for details on the lazy peeling of non-bipartite Boltzmann maps).  Let $ \mathbf{q}= ( q_{k})_{k \geq 1}$  be a non-zero sequence of non-negative real numbers, the so-called \textbf{weight sequence}. We define a $\sigma$-finite measure on all finite bipartite rooted planar maps (with no simplicity restrictions) by the formula 
\begin{eqnarray} \label{eq:defBoltzmann} w_{ \mathbf{q}}( \map) &=& \prod_{ f \in \mathsf{Faces}(\map) \backslash \{ \rootface\}} q_{ \mathrm{deg}(f)/2},  \end{eqnarray}
where $\rootface$ is the root face of the map, i.e.~the face incident to the right  of the root edge, which receives no weight in the above formula.
We define $W^{(\ell)}$ to be the total $w_{\mathbf{q}}$-weight of all bipartite rooted maps with a root face of degree $2 \ell$.
We suppose that $\mathbf{q}$ is \emph{admissible}, i.e.~that $W^{(\ell)}$ is finite for any $\ell\geq 1$. 
In that case we can define the associated Boltzmann distribution $\mathbb{P}^{(\ell)}$ by normalizing $w_{ \mathbf{q}}$ on the set of maps with root face degree $2 \ell$. 
In this context, very general enumeration results for Boltzmann planar maps, see \cite[Lemma 3.13]{CurStFlour} give a "strong ratio limit" theorem, in the sense that there exists an explicit constant $c_{ \mathbf{q}}>1$ such that 
  \begin{eqnarray} \frac{W^{(\ell+1)}}{W^{(\ell)}} \xrightarrow[\ell\to\infty]{} c_{ \mathbf{q}}.   \label{eq:strongratio}\end{eqnarray}
This asymptotic enumerative result enables one to define a half-planar model of the random Boltzmann map with a general boundary, in the sense that $ \mathbb{P}^{(\ell)} \to \mathbb{P}^{(\infty)}$ in distribution for the local topology of maps where $ \mathbb{P}^{(\infty)}$ is a probability measure on the set of half-planar infinite maps, see \cite[Chapter VI]{CurStFlour}. This model of half-planar map is very convenient for the lazy-peeling process since the infinite unexplored region always has law $ \mathbb{P}^{(\infty)}$. Moreover the lazy-peeling process is intimately connected with a random walk on $ \mathbb{Z}$ whose step distribution, denoted by $\nu$, is defined by
 \begin{eqnarray} \nu(k) = \left\{ \begin{array}{ll} q_{k+1} c_{ \mathbf{q}}^{k}  & \mbox{for } k \geq 0 \\
2W^{(-1-k)} c_{ \mathbf{q}}^{k} & \mbox{for } k \leq -1. \end{array} \right.   \label{eq:nu}\end{eqnarray}
A central quantity, popping up e.g.~in the analysis of percolation on $ \mathbb{P}^{(\infty)}$, is the so-called \emph{lazy-gulp} $\gulp$, given by half of the mean number of edges of the boundary swallowed during a peeling step, see \cite[Definition 11]{CurStFlour}. In terms of the measure $\nu$ it reads
$$ \gulp = \sum_{k \geq 1} \nu(-k) (2k-1).$$
In particular, when $\gulp < \infty$, the so-called \emph{dilute phase}, the random half-planar map of law $ \mathbb{P}^{(\infty)}$ has a finite number of pinch points separating the root edge from $\infty$ almost surely. In the \emph{dense phase} $\gulp = \infty$ the half-planar map is made of an infinite tree of finite components separated by pinch points. See \cite[Proposition 6.7]{CurStFlour} and \cite{Richier17b} for the case of stable weight sequences.

\paragraph{Enumeration and local limits of maps with simple boundaries.}
Our first goal in this paper is to achieve the same results as above in the case of maps with a \emph{simple boundary}. Recall that a face $f$ of a (bipartite rooted planar) map $\map$ is \textbf{simple} if it contains no pinch points, i.e.~if we do not visit a given vertex more than once while following the contour of the face $f$. A map $\map$ whose root face is simple will be called $\partial$-simple (not to confuse with simple maps where loops and multiple edges are forbidden). 
Note that the map consisting of a single (root) edge and a single face, is considered a $\partial$-simple map. 

\begin{figure}[!h]
	\begin{center}
		\includegraphics[width=0.2\linewidth]{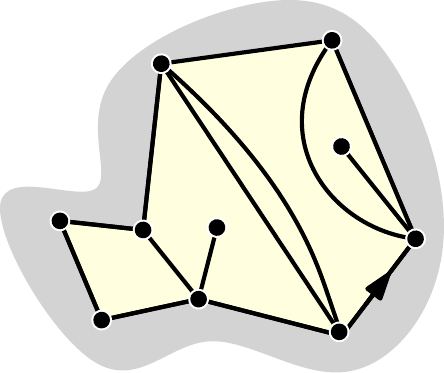}
		\caption{ \label{fig:simplemap} A $\partial$-simple map.}
	\end{center}
\end{figure}
We will add a hat ``$\ \hat{}\ $'' to the above notation when dealing with $\partial$-simple maps. For instance, we denote by $\hat{W}^{(\ell)}$ the total $w_{ \mathbf{q}}$-weight of  the set of all (bipartite rooted planar) map with a simple root face of degree $2 \ell$. We write $ \hat{\mathbb{P}}^{(\ell)}$ for the simple $ \mathbf{q}$-Boltzmann measure obtained by normalizing $w_{ \mathbf{q}}$ on the above set. In order to generalize \eqref{eq:strongratio} to the case of $\partial$-simple maps, we shall suppose that the weight sequence $ \mathbf{q}$ is \emph{critical}, which means here that the $\nu$-random walk oscillates \cite[Proposition~4]{Bud15} or, equivalently, that the variance of the number of vertices of a Boltzmann map under $ \mathbb{P}^{(\ell)}$ is infinite, see \cite[Proposition 4.3]{BCMgasket}.
In this case Proposition \ref{prop:ratiosimple} establishes the existence of $\hat{c}_{ \mathbf{q}}>0$ such that
\begin{equation*} \frac{\hat{W}^{(\ell+1)}}{\hat{W}^{(\ell)}} \xrightarrow[\ell\to\infty]{} \hat{c}_{ \mathbf{q}},
\end{equation*}
which is deduced by purely probabilistic arguments, using the lazy-peeling process, its connection with a random walk and the strong ratio limit theorem. 
Proposition \ref{prop:ratiosimple} is an important ingredient in the construction of the $\partial$-simple analog of $\mathbb{P}^{(\infty)}$:

\begin{theorem}[Local limit of maps with simple boundaries]\label{thm:localsimple}
Let $\qseq$ be a \emph{critical} weight sequence and recall that $\hat{\mathbb{P}}^{(\ell)}$ is the distribution ${\mathbb{P}}^{(\ell)}$ conditioned on maps having a \emph{simple} boundary of perimeter $2 \ell$. 
Then we have the following local weak convergence:	
\begin{equation*} 
\hat{\mathbb{P}}^{(\ell)} \xrightarrow[\ell\to\infty]{(d)}\hat{\mathbb{P}}^{(\infty)},  
\end{equation*}
where $\hat{\mathbb{P}}^{(\infty)}$ is a probability measure supported by infinite maps of the half-plane with a simple boundary.
\end{theorem}
The construction of $\hat{\mathbb{P}}^{(\infty)}$ is easiest in the dilute phase because under $ \mathbb{P}^{(\infty)}$ our half-planar maps have almost surely a unique infinite simple component which is distributed as $\hat{\mathbb{P}}^{(\infty)}$, see \cite{CMboundary} for the quadrangular case. In fact, in the triangular and quadrangular case, this convergence is due to Angel \cite{Ang03} and defines the so-called Uniform Infinite Half-Planar Triangulation or Quadrangulation (with a simple boundary) which have been the subject of numerous studies, see e.g.~\cite{Ang05,AR13,ACpercopeel,CCsaw,baur2016geodesic,Ray13,angel2018half,richier2018incipient,gwynne2017convergence,GM16b,GM16c}. However, there are examples of critical weight sequences $ \mathbf{q}$ for which $\mathbb{P}^{(\infty)}$ has no infinite simple component (see \cite{Richier17b}) and so $\hat{\mathbb{P}}^{(\infty)}$ cannot be defined using a pruning procedure and we shall rather rely on a delicate Doob $h$-transformation. 
We then study the simple peeling process under $ \hat{\mathbb{P}}^{(\infty)}$ and compare it with the lazy-peeling process under $ \mathbb{P}^{(\infty)}$. One particular finding, is the fact that the associated simple gulp $ \sgulp$, defined as half of the mean number of edges of the boundary swallowed during a simple-peeling step, is identical to that of the lazy-peeling process! This surprising fact is proved by analyzing Bernoulli percolations but begs for a more straightforward explanation.
 
\paragraph{Application to percolation thresholds.} Since the pioneering paper of Angel \cite{Ang03}, there has been a lot of work on Bernoulli percolations on random maps (including among others \cite{Ang05,ACpercopeel,Ray13,MN13,richier2018incipient,gwynne2017convergence,holden2019convergence,bernardi2018percolation,CKperco,BCMgasket,bjornberg2015site,SheHC}), and the peeling process turned out to be the tool of choice in their analysis. Recall that Bernoulli site/bond/face-percolation consists in coloring independently the vertices/edges/faces of a map in black with probability $p\in(0,1)$ and white otherwise. Using the \emph{simple peeling}, the critical percolation thresholds (for existence of an infinite black cluster) have been computed in \cite{Ang03,ACpercopeel,Ric15} in the case of half-plane triangulations or quadrangulations with a simple boundary. An alternative computation is provided in \cite{CurStFlour,BCMcauchy,curien2018duality} using the \emph{lazy-peeling} but only for bond and face percolations. In this work, we complete the picture by determining the explicit thresholds for all bipartite Boltzmann maps:

\begin{theorem}[Percolation thresholds]\label{thm:siteperco}\noindent Let $\qseq$ be a critical weight sequence in the dilute phase $\gulp < \infty$ and $\nu$ as defined in \eqref{eq:nu}. Then under $\hat{\mathbb{P}}^{(\infty)}$ or ${\mathbb{P}}^{(\infty)}$ the critical percolation threshold for Bernoulli site/bond/face-percolation is almost surely constant and equal to 
\begin{equation}
p_{c, \mathrm{site}}= 1 - \frac{(\sum_{k=1}^\infty \nu(-k) )^2}{ 2\nu(-1)  \gulp }, \quad p_{c, \mathrm{bond}}= 1 - \frac{1}{ \gulp +1},\quad  p_{c, \mathrm{face}}= \frac{1}{2}\left(1 + \frac{1}{ 2 \gulp +1}\right).
\end{equation}
Furthermore, there is no percolation at criticality.
\end{theorem}
The truly original part of this theorem is the computation of $p_{c, \mathrm{site}}$. Indeed, the universal forms for $ p_{c, \mathrm{bond}}$ and $p_{c, \mathrm{face}}$ can already be found in \cite{ACpercopeel} (with the notation $\delta^{*}= 2 \sgulp$) and in \cite{BCMcauchy,curien2018duality}. The site percolation threshold was computed by Richier \cite{Ric15} in the case of quadrangulations but the formula involved another quantity besides the $\gulp$ related to the simple peeling (see also \cite{bjornberg2015site}). Our finding is that this quantity can be explicitly expressed in terms of the lazy-peeling and gives rise to the formula above. As an example we record in Table \ref{tbl:values} the thresholds for $2p$-angulations, and refer to \cite[Section 8]{Bud15} and \cite[Section~3.4]{CurStFlour}  for explicit computations. Several surprising identities show up along the way, for instance in Proposition~\ref{prop:simplegulp} and Lemma~\ref{lem:gulpcalc}, and we hope that these will find a better combinatorial explanation in the future.

\begin{table}
\centering
\begin{tabular}{l | c c c c}
 & quadrangulations & $6$-angulations & $8$-angulations & $2p$-angulations \\\hline
$p_{c, \mathrm{site}}$ & $\frac{5}{9}$ & $\frac{76}{125}$ & $\frac{5197}{8085}$ & $1-\frac{2p}{p-1} \frac{\left(\frac{\expo}{1-2p}+1\right)^2}{\expo-1}$ \\
$p_{c, \mathrm{bond}}$ & $\frac{1}{3}$ & $\frac{5}{11}$ & $\frac{11}{21}$ & $\frac{\expo-1}{\expo+1}$ \\
$p_{c, \mathrm{face}}$ & $\frac{3}{4}$ & $\frac{11}{16}$ & $\frac{21}{32}$ & $ \frac{1}{2}(1 + \frac{1}{\expo})$
\end{tabular}
\caption{\label{tbl:values}Explicit values of the percolation thresholds for $2p$-angulations. The \emph{exposure} $\expo = 2 \gulp +1$ in the last column is given by $\expo = 4^{p-1} / \binom{2p-2}{p-1}$. }
\end{table}

\medskip 

\subsection*{Acknowledgements} 
We thank the anonymous referees for useful suggestions and corrections.
The first author acknowledges support of the START-UP 2018 programme with project number 740.018.017, which is financed by the Dutch Research Council (NWO).

\section{Peeling explorations}
In this section we briefly recall the basics about explorations of maps. They rely on two notions of a submap, which is different in the case of the lazy and the simple peeling processes. The reader is referred to \cite{CurStFlour} for more details, especially in the case of the lazy-peeling process (see also the expositions in \cite{Bud15,BC16}). If $\map$ is a planar map, we denote by $|\partial \map|$ the half-perimeter of its root face (this is an integer since we are working with bipartite maps).

\subsection{Submaps} Let $\map$ be a (bipartite rooted planar) finite or infinite map and $\map^\dagger$ the corresponding dual map.
The vertex of $\map^\dagger$ dual to the root face of $\map$, i.e. the face incident to the right of its root edge, is called the \emph{origin} of $\map^\dagger$. 
Let $ \mathfrak{e}^\circ$ be a finite connected subset of edges of $\map^\dagger$ that is incident to the origin of $\map^\dagger$ (the letter ``$\mathfrak{e}$'' stands for \textbf{explored}). We associate with $\mathfrak{e}^\circ$ a planar map $\mathfrak{e}$ which, roughly speaking, is obtained by gluing the faces of $\map$ that correspond to the (dual) vertices incident to $\mathfrak{e}^\circ$ (or just the root face if $\mathfrak{e}^\circ$ is empty) along the (dual) edges of $ \mathfrak{e}^\circ$. 
Note that several vertices of $\mathfrak{e}$ may correspond to the same vertex in $\map$.
We associate a second map $\hat{ \mathfrak{e}}$ to $\mathfrak{e}^\circ$ by identifying such vertices in $\mathfrak{e}$.
In other words, $\hat{ \mathfrak{e}}$ is obtained from $\map$ by retaining the faces and vertices incident to $\mathfrak{e}^\circ$ but cutting along the edges of ${\map}$ whose dual edges are not in $ \mathfrak{e}^{\circ}$. See Figure \ref{fig:simple-submap} below. 

\begin{figure}[!h]
	\begin{center}
		\includegraphics[width=.85\linewidth]{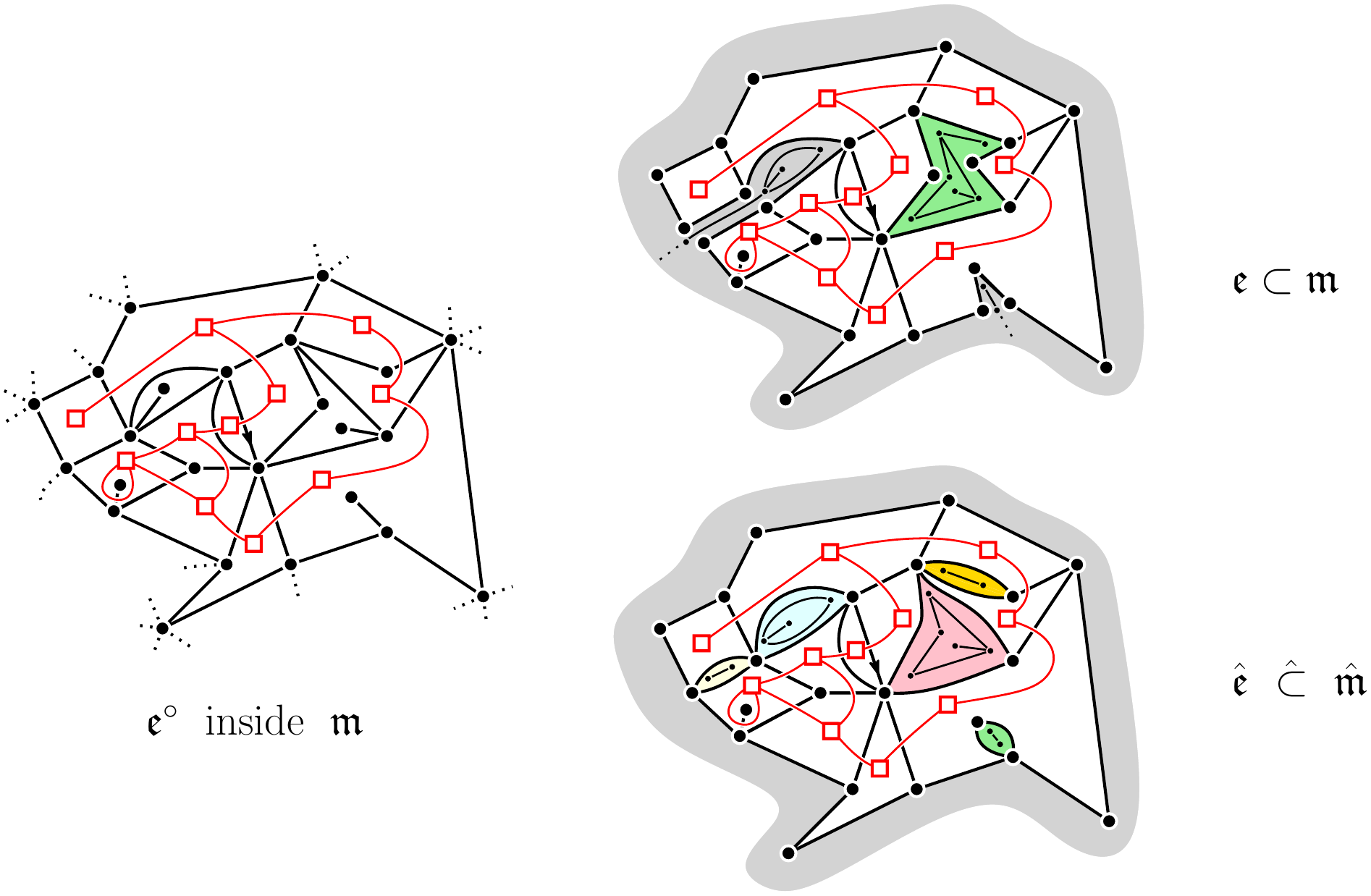}
		\caption{ \label{fig:simple-submap} Illustration of the two notions of submaps. On the left, part of a map ${\map}$ with a simple boundary of perimeter $2$ is shown in black together with a connected subset $\mathfrak{e}^{\circ}$ of its dual in red. 
		On the right, lazy-submaps (top) and simple-submaps (bottom) $\mathfrak{e}$ and $\hat{\mathfrak{e}}$ obtained by cutting along all edges not in $\mathfrak{e}^{\circ}$ by keeping or not the identification of vertices inside $ \hat{\map}$. Each shaded region corresponds to a different hole, for which we have indicated the map that is to be glued in the hole in order to recover $\map$. }
	\end{center}
\end{figure}

Both $ \mathfrak{e}$ and $ \hat{ \mathfrak{e}}$ are planar maps with several distinguished faces, called the \emph{holes}, that correspond to the connected components of $\map$ that are cut out when building $ \mathfrak{e}$, resp.~$ \hat{ \mathfrak{e}}$. Notice that the holes are simple and, in the case of $ \mathfrak{e}$, they do not share any vertices. In the case of $ \hat{ \mathfrak{e}}$, the holes can share vertices, but those vertices together cannot disconnect the map, because $\mathfrak{e}^\circ$ is connected. We say that $ \mathfrak{e}$ (resp.~$ \hat{ \mathfrak{e}}$) is a \emph{lazy-submap} (resp.~\emph{simple-submap}) of $\map$ and write  $$ \mathfrak{e} \subset \map \qquad (\mbox{resp.~} \quad  \hat{\mathfrak{e}}  \hatsubset \map),$$ since $\map$ can be recovered from $ \mathfrak{e}$ (resp.~ $ \hat{ \mathfrak{e}}$) by gluing\footnote{To perform this gluing operation, we implicitly assume that an oriented edge is distinguished on the boundary of each hole and that the root edges of the maps are identified. We will not mention this further, since these edges can be arbitrarily chosen using a deterministic procedure. The gluing operation is rigid, in the sense that the maps filling-in the holes of $ \hat{ \mathfrak{e}}$, once rooted, are uniquely determined.} inside each hole a map with a general boundary (resp.~simple boundary). 
Observe that in the case of a simple-submap under this gluing operation two edges on the boundary of the same hole are never identified, unless the hole is of degree two and one glues the $\partial$-simple map consisting of a single edge into the hole.
The union of the boundaries of the holes is called the \emph{active boundary} of $ \mathfrak{e}$ (resp. $\hat{ \mathfrak{e}}$). In the rest of the paper, for any submap  $ e$, lazy or simple, its Boltzmann weight is always
$$  w_{ \mathbf{q}}( e) = \prod_{ f \in \mathsf{Faces}(e) \backslash \{ \rootface, \text{holes}\}} q_{ \mathrm{deg}(f)/2}.$$

\begin{remark}[Root transform] In the case of the simple-submap we will always implicitly assume that $\map$ is $\partial$-simple. In particular, if $\mathfrak{e}^{\circ}$ is empty then the associated submaps $\mathfrak{e}_{0}$ and $\hat{\mathfrak{e}}_{0}$ both consist of a simple face (corresponding to the root face of $\map$) and a hole of the same degree glued together. This is not a demanding assumption, since any bipartite map $\map$ can be seen as a $\partial$-simple map of the $2$-gon after splitting the root edge. See \cite[Figure 3.2, page 43]{CurStFlour}.
\end{remark}

\subsection{Explorations}
We can now define what we mean by a peeling exploration of a map. Formally, this depends on a function $ \mathcal{A}$, called the \textbf{peeling algorithm}, which associates to any planar map $\mathfrak{e}$ with holes an edge of the active boundary of $\mathfrak{e}$ or the symbol \Cross, which we interpret as the instruction to end the exploration. 
The lazy-peeling of the map $\map$ with algorithm $\mathcal{A}$ is the sequence of lazy-submaps $$ \mathfrak{e}_{0} \subset \mathfrak{e}_{1} \subset \cdots \subset \mathfrak{e}_{n} \subset \cdots \subset \map$$  associated to the sequence $(\mathfrak{e}_{i}^\circ)_{i \geq 0}$ of growing connected subsets of edges of $\map^\dagger$, where $\mathfrak{e}_{0}^{\circ}$ is the root face and such that $\mathfrak{e}_{i+1}^\circ$ is obtained from $\mathfrak{e}_{i}^\circ$ by adding the edge dual to $\mathcal{A}(\mathfrak{e}_{i})$ inside $\map$ (unless the exploration has stopped). Replacing lazy-submaps by simple-submaps, we can analogously define the simple peeling exploration of $\map$ with algorithm $\mathcal{A}$ as the growing sequence of simple-submaps
	$$ \hat{\mathfrak{e}}_{0}  \hatsubset   \hat{\mathfrak{e}}_{1} \hatsubset\cdots  \hatsubset   \hat{\mathfrak{e}}_{n} \hatsubset  \cdots \hatsubset    \map,$$ associated to the sequence $(  \mathfrak{e}_{i}^{\circ})_{i \geq 0}$ where $  \mathfrak{e}_{0}^{\circ}$ is the root face and $ \mathfrak{e}^{\circ}_{n+1}$ is obtained by adding the edge dual to $ \mathcal{A}( \hat{ \mathfrak{e}}_{n})$ inside $\map$. 
	
	The move $ \mathfrak{e}_{n} \to \mathfrak{e}_{{n+1}}$ which is obtained by \emph{peeling the edge $ \mathcal{A}( \mathfrak{e}_{n})$} in the lazy-exploration can be classified into two topologically different types (provided that $  \mathcal{A}( \mathfrak{e}_{n})\ne\mbox{ \Cross}$), depending on  the face $ \mathsf{F}$ of $\map$ adjacent to $ \mathcal{A}( \mathfrak{e}_{n})$ and located on the other side of $  \mathcal{A}( \mathfrak{e}_{n})$ with respect to $ \mathfrak{e}_{n}$:
	\begin{itemize}
\item \emph{Event $ \mathsf{C}_{k}$:}  the face $\mathsf{F}$ is not a face of $\mathfrak{e}_{n}$ and has degree $2k$. Then  $ \mathfrak{e}_{n+1}$ is obtained by gluing $\mathsf{F}$ to $\mathcal{A}( \mathfrak{e}_{n})$ without performing the possible identifications of its other edges inside $\map$.

\item  \textit{Event $\mathsf{G}_{k_{1},k_{2}}$:} the face $\mathsf{F}$ is already a face of $\mathfrak{e}_{n}$. In this case, the edge $ \mathcal{A}( \mathfrak{e}_{n})$ is identified in $\map$ with another edge $a'$ on the boundary of the same hole where $k_{1}$ (resp.~$k_{2}$) is half of the number of edges on the boundary of the hole strictly between $ \mathcal{A}( \mathfrak{e}_{n})$ and $a'$ when turning in clockwise order around the hole, and $ \mathfrak{e}_{n+1}$ is the map after this identification in $\mathfrak{e}_n$.
\end{itemize}

Note that when $k_{1}>0$ and $k_{2}>0$ the event $\mathsf{G}_{k_{1},k_{2}}$ results in the splitting of a hole into two holes. If $k_{1}=0$ or $k_{2}=0$ the corresponding hole of perimeter $0$ is simply a vertex of the map. In particular, the event $\mathsf{G}_{0,0}$ results in the disappearance of a hole (of degree $2$).

\begin{figure}[!h]
 \begin{center}
 \includegraphics[height=8cm]{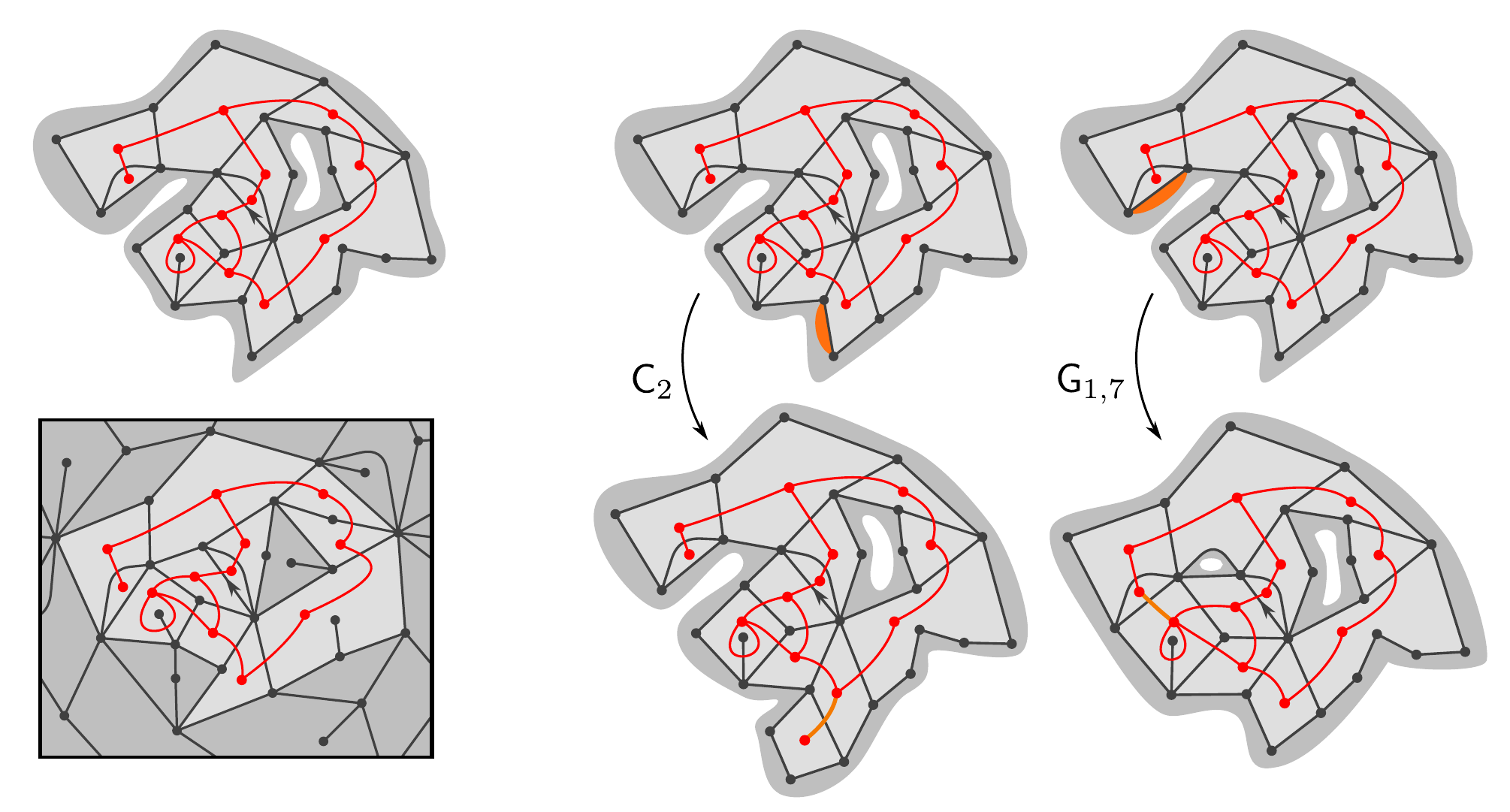} \caption{\label{fig:peel}Illustration of the different lazy-peeling events. The left column depicts the submap $\mathfrak{e}\subset \map$ (top) and the map $\map$ with $\mathfrak{e}^\circ$ superimposed in red. The center and right columns represent two different peeling events ($\mathsf{C}_{2}$ and $\mathsf{G}_{1,7}$) depending on the edge to be peeled (in thick orange). }
 \end{center}
 \end{figure}

\begin{figure}[!h]
	\begin{center}
		\includegraphics[width=0.8\linewidth]{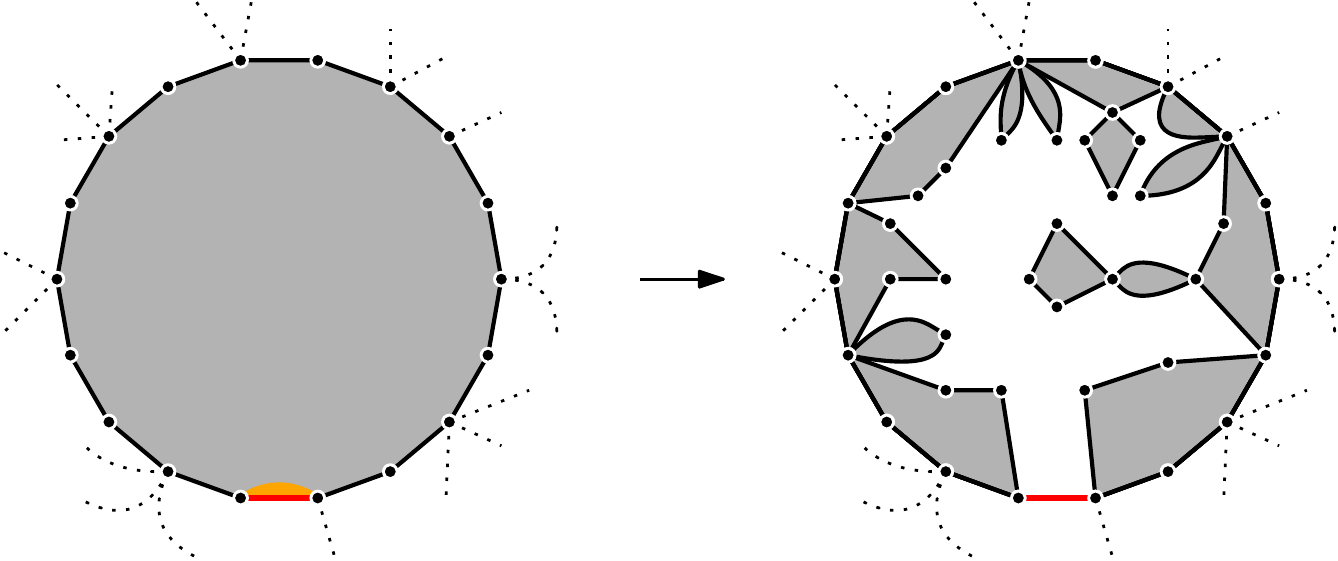}
		\caption{ \label{fig:transition-simple-peeling}An example of move from $ \hat{\mathfrak{e}}_{n}$ to $ \hat{\mathfrak{e}}_{n+1}$. The edge to peel is in orange and the hole of perimeter $18$ (in gray) is split into $14$ new holes of perimeters $6,2,6,6,2,2,4,4,2,2,6,2,4,6$ after discovering a new face (in white) of perimeter $38$.}
	\end{center}
\end{figure}

For the simple-peeling, however, the moves $ \hat{ \mathfrak{e}}_{n} \to \hat{ \mathfrak{e}}_{n+1}$ cannot be classified into finitely many types because many different topological situations can occur, like the one in Figure \ref{fig:transition-simple-peeling}. We shall not try to formalize these possible transitions for the time being.

\medskip 

\emph{Filled-in explorations.} In this paper, it is convenient to consider explorations where the submaps maintain a \emph{unique} hole at each step. Those explorations, denoted by $ (\overline{ \mathfrak{e}}_{n} : n \geq 0)$ in the lazy case and $(\hatbar{\mathfrak{e}}_{n} : n \geq 0)$ in the simple case are described as above, except that when the peeling of an edge in $ \hatbar{\mathfrak{e}}_{n}$ or in $\hatbar{\mathfrak{e}}_{n}$ results in the splitting of a hole into two or more holes (i.e.~for event $ \mathsf{G}_{k_{1},k_{2}}$ with $k_{1},k_{2} \geq 1$ in the lazy case), then we fill in all but one (specified by the algorithm) with their respective parts of $\map$. See Section 4.1.3 of \cite{CurStFlour} for more details. 

\subsection{Peeling of Boltzmann maps}

When applied to a random Boltzmann map, the above peeling explorations turn out to be Markov chains whose probability transitions can be explicitly described in terms of $W^{(\ell)}$ and $ \hat{W}^{(\ell)}$. More precisely, if $( \mathfrak{e}_{n} : n \geq 0)$ is a lazy-peeling exploration of the map $\map$ distributed according to $ \mathbb{P}^{(\ell)}$, then for each $n \geq 0$, conditionally on $ \mathfrak{e}_{n}$, the maps filling in the holes of $ \mathfrak{e}_{n}$ inside $\map$ are independent $\mathbf{q}$-Boltzmann maps with the appropriate perimeters. In particular, this enables us to compute the probability of the events $ \mathsf{C}_{k}$ and $ \mathsf{G}_{k_{1},k_{2}}$ above: conditionally on the past exploration, if  $ \mathcal{L}_{n}$ is the half-perimeter of the hole on which $ \mathcal{A}( \mathfrak{e}_{n})$ is selected, then the events $ \mathsf{C}_{k}$ and $\mathsf{G}_{k_{1},k_{2}}$ (where $k \geq 1$ and $k_{1}+k_{2}+1=\mathcal{L}_{n}$ with $k_{1}, k_{2} \geq 0$) occur respectively with probabilities
$$q_{k} \frac{W^{(\mathcal{L}_{n}+k-1)}}{W^{(\mathcal{L}_{n})}} \qquad \mbox{ and } \qquad \frac{W^{(k_{1})}W^{(k_{2})}}{W^{( \mathcal{L}_{n})}}.$$
We refer to \cite[Section 4.2]{CurStFlour} for details.
With the help of \eqref{eq:strongratio}, we note that the limit $\mathcal{L}_{n} \to \infty$, while keeping $k_1$ or $k_2$ fixed for the second probability, gives rise to the measure $\nu$ defined in \eqref{eq:nu}. 

An analogous statement holds for the simple-peeling exploration after adding the appropriate hats ``$\hat{\ \ }$''. For example, the transitions of the form of Figure \ref{fig:transition-simple-peeling} happen with probability 
$$ \frac{1}{ \hat{ {W}}^{(\ell)}}q_{k} \prod_{i} \hat{W}^{(\ell_{i})},$$
where $k$ is the half-perimeter of the face we reveal, $\ell$ the half-perimeter of the hole on which the edge to peel has been selected, and $\ell_{i}$ are the half-perimeters of the holes created. If $\ell=1$, we might not reveal any face and glue the two sides of the $2$-gon together and close it.  This happens with probability $1 / \hat{W}^{(1)}$. Here also, conditionally on the past exploration, the maps filling in the holes inside $\map$ are independent and distributed as simple Boltzmann maps with the correct perimeter.  \medskip 

In particular, when considering filled-in explorations of Boltzmann maps, the holes we fill in are distributed as Boltzmann maps of law $\mathbb{P}^{(\ell)}$ or $\hat{ \mathbb{P}}^{(\ell)} $ with the correct perimeters. 

\section{Weight of maps with simple boundaries}
\label{sec:simpleenu}
In this section we gather the basics on $\partial$-\emph{simple} maps and in particular prove the asymptotic estimate of Proposition \ref{prop:ratiosimple} which is the analog of  \eqref{eq:strongratio} in the general boundary case. This will serve as a key ingredient in the proof of Theorem \ref{thm:localsimple}.

\subsection{Core decomposition}

If $\map$ is a map with a general boundary, one can consider the map with a simple boundary obtained by keeping the simple boundary component carrying the root edge, see Figure \ref{fig:core}. The simple component carrying the root edge will be called the \emph{core} and denoted by $ \mathsf{Core}(\map)$. This well-known decomposition \cite{BIPZ78,BG09,ABM16} yields a way to relate $\hat{W}^{(\ell)}$ to their non-simple analogs, see \cite[Lemma 2.6]{Richier17b}:
\begin{figure}[!h]
	\begin{center}
		\includegraphics[width=14cm]{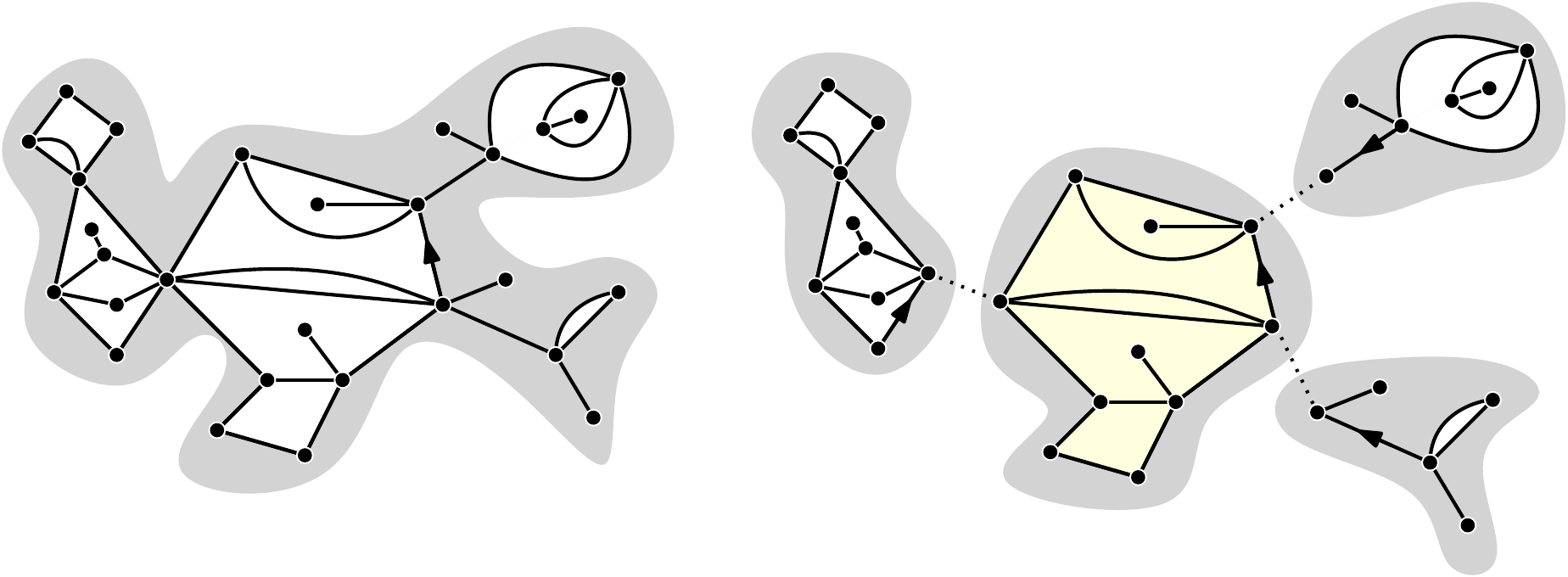}
		\caption{ \label{fig:core}The decomposition of a map with a general boundary. On the left, a map $\map$ with a general boundary. On the right, we see the simple component of the root edge (in yellow) together with the maps with general boundaries attached to it.}
	\end{center}
\end{figure}

\begin{proposition}\label{prop:simplegenfun}  Let $ \qseq$ be a weight sequence. If $ \mathrm{W}(z) = \sum_{\ell \geq 0} W^{(\ell)}z^{\ell}$ and $\hat{ \mathrm{W}}(z) = \sum_{ \ell \geq 0}  \hat{W}^{(\ell)}z^{\ell}$ then we have the following equality between (formal)  power series:
	$$ \hat{ \mathrm{W}}\left(z ( \mathrm{W}(z))^{2}\right) =  \mathrm{W}(z).$$
\end{proposition}

The strong ratio limit \eqref{eq:strongratio} implies that the radius of convergence of $\mathrm{W}$ is $ \frac{1}{c_{ \qseq}}$ and therefore $\mathrm{W}(1/c_{ \qseq}) > 1$. Although we will not use this, \cite[Lemma 5.2]{CurStFlour} implies the upper bound $\mathrm{W}(1/c_{ \qseq})\leq\frac{c_{\mathbf{q}}}{2}$. It follows from the previous proposition that the radius of convergence of $ \hat{ \mathrm{W}}$ is \emph{at least} $ \frac{1}{ \hat{c}_{ \qseq}}$ where 
\begin{eqnarray} \label{def:cqhat} \hat{c}_{ \qseq} = \frac{c_{ \qseq}}{ \mathrm{W}^{2}( \frac{1}{c_{ \qseq}})},  \end{eqnarray}
and that $ \hat{ \mathrm{W}}( {1}/{ \hat{c}_{ \qseq}}) = \mathrm{W}(1/c_{ \qseq})$, see also \cite{Richier17b}. In the following we write $ \mathrm{W}_{c} = \mathrm{W}(1/ c_{ \qseq})$ and $\mathrm{W}'_{ \mathrm{c}} = \mathrm{W}'(1/ c_{ \qseq})$ to simplify notation. We shall prove in this section the following key asymptotic estimate,  which will eventually follow from the strong ratio limit theorem for random walks:

\begin{proposition}[Strong ratio limit for $\hat{W}^{(\ell)}$] \label{prop:ratiosimple}Let $ \qseq$ be an admissible \emph{critical} weight sequence. Then we have $$ \displaystyle \lim_{\ell \to \infty} \frac{ \hat{W}^{(\ell+1)}}{ \hat{ W}^{(\ell)}} = \hat{c}_{ \qseq}.$$
\end{proposition}
\begin{remark} In \cite{Richier17b}, Richier provides related estimates on the function $ \hat{ \mathrm{W}}$ using Tauberian theorems in the case when $ \mathbf{q}$ is a regular varying weight sequence. Those estimates imply in some cases that the radius of convergence of $ \hat{ \mathrm{W}}$ is equal to $\frac{1}{ \hat{c}_{ \qseq}}$. 
\end{remark}

\begin{remark}\label{rem:subcritical} For subcritical weight sequences we may have  $\lim_{\ell \to \infty} \frac{ \hat{W}^{(\ell+1)}}{ \hat{ W}^{(\ell)}} > \hat{c}_{ \qseq}$ and this is related to the fact that the perimeter of the core of a Boltzmann map of law $\mathbb{P}^{(\ell)}$ may have an exponential tail as $\ell \to \infty$. We shall thus restrict to critical sequences in what follows. See \cite[Remark 2.8]{Richier17b} for the case of quadrangulations.\end{remark}

\begin{remark} Notice that when $ \qseq$ is critical then it is easy to define  infinite $ \qseq$-Boltzmann map of the \emph{plane} with simple boundary of perimeter $2 \ell$ by conditioning the infinite Boltzmann maps of the plane with perimeter $2 \ell$ (see \cite[Chapter VII]{CurStFlour}) to have a simple boundary (this is an event of positive probability). The difference between maps with simple boundary and maps with general boundary will thus manifest itself more dramatically when dealing with half-planar limits.\end{remark} 

\hrulefill \textit{ In the rest of the paper we shall assume that $ \qseq$ is critical.} \hrulefill \\ 

The rest of this section is devoted to proving Proposition \ref{prop:ratiosimple}, and for this we need to develop a peeling algorithm to explore the core of a planar map whose perimeter is random.

\subsection{Free Boltzmann map and exploration of the core}
\label{sec:freeboltz}
In the following it will be convenient to work with Boltzmann maps with a general boundary whose perimeter is not fixed.
We introduce the probability measure
\begin{equation}\label{eq:Pfreedef}
\mathbb{P}^{( \mathrm{free})} =  \frac{1}{ \mathrm{W}_{c}} \sum_{\ell \geq 0} W^{(\ell)} c_{ \qseq}^{-\ell} \cdot \mathbb{P}^{(\ell)}.
\end{equation}
The random planar map of law $ \mathbb{P}^{( \mathrm{free})}$ thus has a random boundary length, but, conditionally on this length, is a standard Boltzmann map (with general boundary).
We allow the free Boltzmann map to have boundary length $0$, which occurs with probability $1/\mathrm{W}_{c}$, corresponding to the unique map $\dagger$ consisting of a single vertex (and no edges).
From the decomposition of Figure \ref{fig:core} or Proposition \ref{prop:simplegenfun} we deduce that
\begin{align} 
\mathbb{P}^{( \mathrm{free})}(|\partial   \mathsf{Core}( \map)| = \ell) &= \frac{1}{W(1/c_{\mathbf{q}})} [y^\ell] \hat{W}\left(y\,\frac{W(1/c_{\mathbf{q}})^2}{c_{\mathbf{q}}}\right)  \nonumber\\
&=  \hat{W}^{(\ell)} c_{\qseq}^{-\ell} (W_{c})^{2\ell-1} =  \frac{\hat{W}^{(\ell)} \hat{c}_{\qseq}^{-\ell}}{W_{c}}. \label{eq:coreperimeter}  \end{align}
To evaluate the probability on the left-hand side, we design a filled-in peeling algorithm that explores the boundary of the core of a map and computes its perimeter. 
Contrary to the lazy peeling of \cite{Bud15}, the perimeter of the maps we consider is itself random and should be thought of as unknown during the peeling. To cope with this difficulty and also deal with maps of the half-plane, we slightly modify our notion of submap (we keep the same name but change the notation).

\begin{figure}[!h]
 \begin{center}
 \includegraphics[width=14cm]{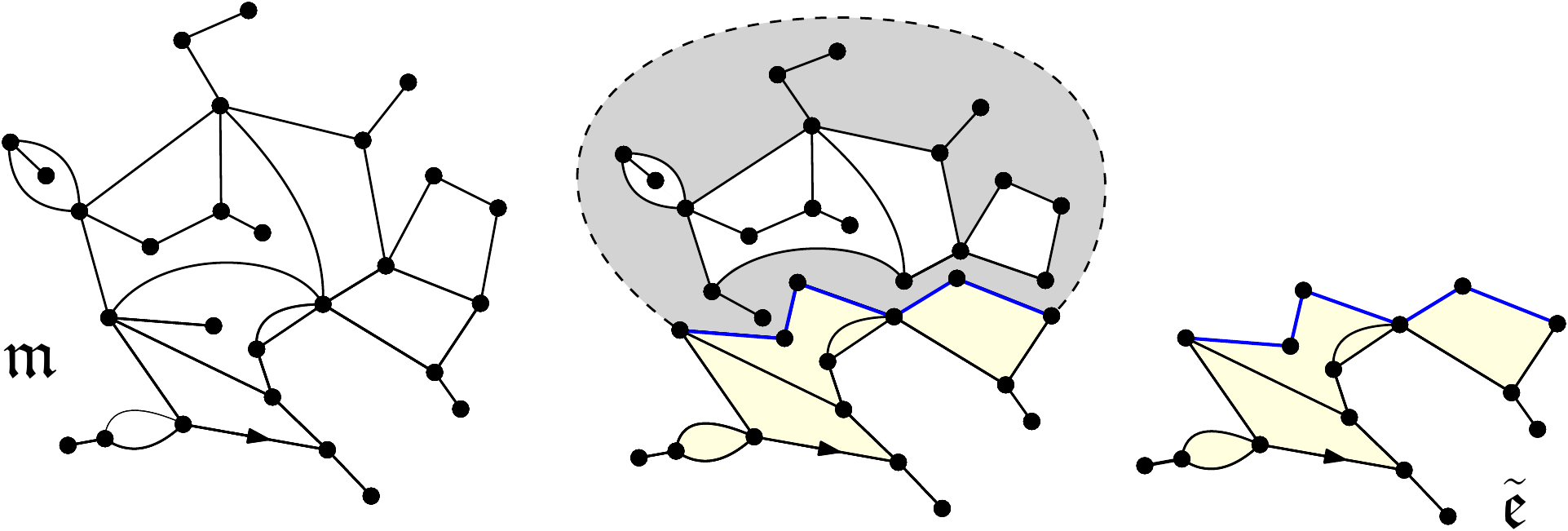}
 \caption{  \label{fig:submaprandom}Illustration of the notion of submap (with one hole) of a map  with an unknown perimeter. On the left, a map $ \mathfrak{m}$ with a general boundary which can be obtained (center picture) as the gluing of a map with general boundary inside the unique hole of  $ \tilde{ \mathfrak{e}}$.}
 \end{center}
 \end{figure}

A submap $ \tilde{ \mathfrak{e}}$ with one hole, will be a map whose boundary is made of two segments (possibly sharing vertices and edges): one with no simplicity condition (shown in black on the right of Figure \ref{fig:submaprandom}) and containing the root edge called the \emph{internal boundary}, and another one (in blue on the right of Figure \ref{fig:submaprandom}) which must be simple and called the \emph{exposed boundary}. We imagine that the unique hole of $ \tilde{ \mathfrak{e}}$ is obtained by linking the extremities of the exposed boundary by a dotted line, see Figure \ref{fig:submaprandom}. We say that $ \tilde{ \mathfrak{e}}\Subset \mathfrak{m}$ if $ \mathfrak{m}$ can be obtained from $ \tilde{ \mathfrak{e}}$ by gluing a map with a general boundary in its hole. The only constraint on the glued map is that it has a perimeter larger than or equal to the length of the exposed boundary (see Figure \ref{fig:submaprandom}). The exposed boundary of $ \tilde{ \mathfrak{e}}$ is denoted by $\partial^* \tilde{ \mathfrak{e}}$ and its internal boundary by $ \partial \tilde{ \mathfrak{e}}$.

In this context, the spatial Markov property of the measure $ \mathbb{P}^{( \mathrm{free})}$ can be stated as follows: if $ \tilde{ \mathfrak{e}}$ is a submap with one hole as above, then conditionally on $\tilde{ \mathfrak{e}} \Subset  \mathfrak{m}$,  the remaining map $ \mathfrak{m} \backslash \tilde{ \mathfrak{e}}$ (i.e.\ the unambiguous map glued on $\tilde{ \mathfrak{e}}$ to obtain $ \mathfrak{m}$, properly rooted) is distributed as $ \mathbb{P}^{( \mathrm{free})}$ conditioned on having a perimeter at least equal to the length of its exposed boundary. We can now describe the filled-in peeling algorithm used to explore the core of a map with a general boundary and undetermined perimeter:

\begin{figure}
	\begin{center}
		\includegraphics[width=.9\linewidth]{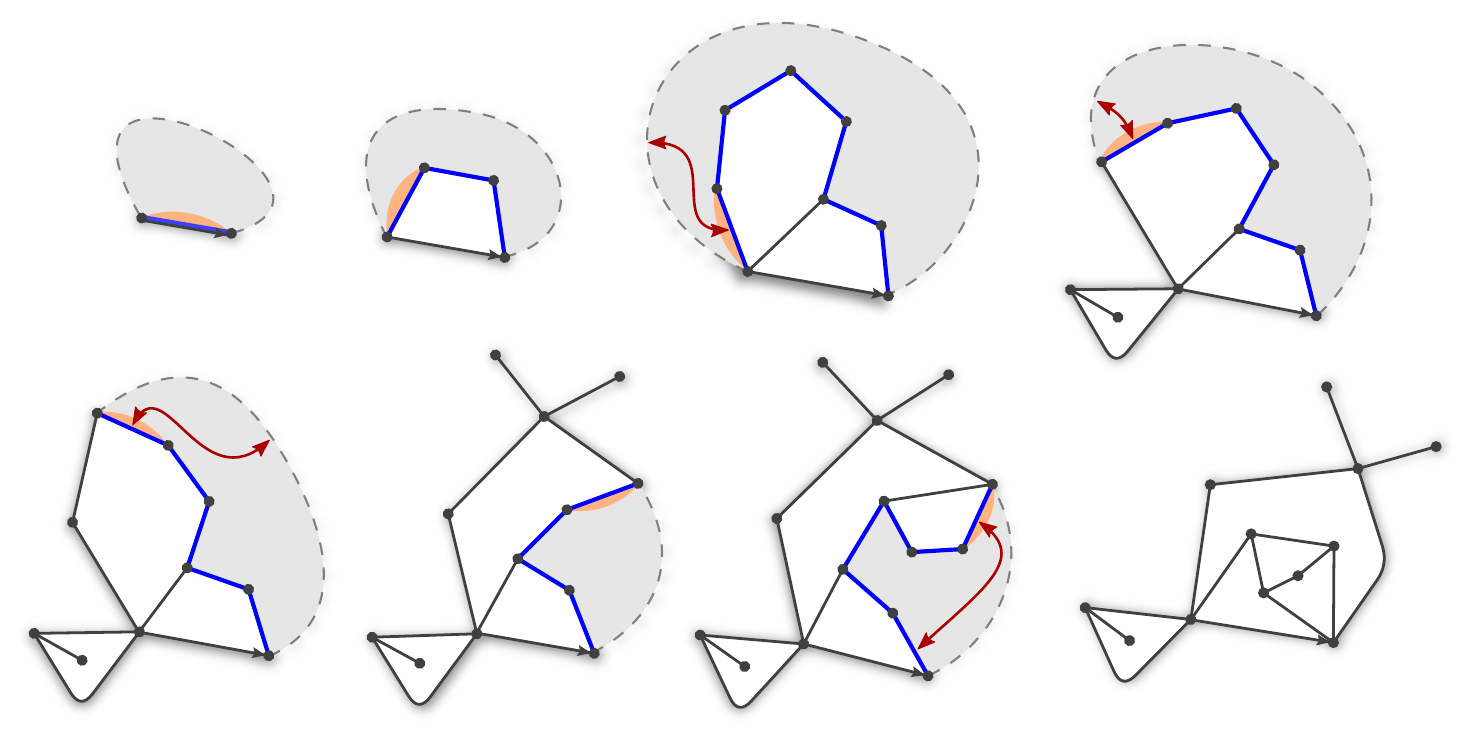}
	\end{center}
	\caption{\label{fig:freepeeling} Example of application of Algorithm $\mathcal{A}_{\text{core}}$ to the map on the bottom right. In the figures, the dashed line represents a sequence of zero or more edges. The exposed boundary is shown in blue, while the next edge to peel is indicated in orange. The arrows indicate which edges are glued in the events of type $\mathsf{G}_{k_1,k_2}$.}
\end{figure}

\begin{center}
	\fbox{\begin{minipage}{15cm}
			\paragraph{Algorithm $\mathcal{A}_{\text{core}}$:} 
			Start with $\tilde{\mathfrak{e}}_{0}$ to be a single edge, seen as having an internal and exposed boundary of length $1$. At each step $n \geq 0$, we peel at the left-most edge of the exposed boundary of $\tilde{ \mathfrak{e}}_n$ and we fill-in all the holes that are not incident to the endpoint of the root edge. The exploration stops when the length of the exposed boundary drops to $0$. This happens either when the exposed boundary is of length $1$ and the peeled edge is identified with another edge of the hole, or when the peeled edge is identified with the right-most edge of the exposed boundary (see the last step in Figure \ref{fig:freepeeling}).
			By convention, if $\map=\dagger$ is the vertex map, the exploration stops immediately.
		\end{minipage}}
	\end{center}

Let us write $(E_{n})_{n \geq 0}$ for the length of the exposed boundary during the exploration with algorithm $ \mathcal{A}_{ \mathrm{core}}$ so that $ E_{0}=1$ (if $\map$ is not the vertex map) and the exploration stops at the first time $\tau$ when $E_{\tau}=0$. We denote by $D$ the number of $-1$ steps that the process $(E)$ performed until time $\tau$.
Note that by construction the left-most edge on the exposed boundary, the peeled edge, always belongs to the core of $\map$.
Each $-1$ step corresponds to an event in which this peeled edge is glued to an edge that is not on the exposed boundary, which therefore must correspond to an edge on the boundary $\partial \mathrm{Core}(\map)$ of the core.
Vice versa, each edge of $\partial \mathrm{Core}(\map)$ except for the root edge must be encountered in the exploration, implying that
	\begin{eqnarray} \label{eq:perimetercore} \mbox{ the perimeter $ 2|\partial \mathrm{Core}(\map)|$ is equal to $D+1$}. \end{eqnarray}
	
\subsection{Evolution of the exposed boundary}	In the next proposition we will give the law of the process $(E_{n}: n \geq 0)$ under $ \mathbb{P}^{( \mathrm{free})}$. But before this, let us introduce some notation and context. In particular, we remind the reader of the classic $h$-transformation of Doob.

	\paragraph{$h$-transformation.} Suppose that  $ p(x,y)$ are probability transitions of a discrete Markov chain on a countable state space $ \Omega$ and that  $h : \Omega \to \mathbb{R}_+$ is a non-negative function which is harmonic and positive on $ A \subset \Omega$, i.e.
$$ h(x) >0, \quad \mbox{ and } \quad h(x) = \sum_{y \in \Omega} p(x,y) h(y), \qquad \forall x \in A,$$
and that $h$ is zero on $\Omega \backslash A$. Under these circumstances, one can define a new transition kernel $q$ on  $A$ by the formula:
 \begin{eqnarray} \label{eq:htransform} q(x,y) = \frac{h(y)}{h(x)} p(x,y), \quad x\in A, y \in \Omega.  \end{eqnarray}
It is plain from the harmonicity of $h$ on $A$ that $ q$ indeed defines a transition kernel. The Markov chain obtained by starting in $A$ is called the \emph{Doob $h$-transform of $p$}. Since $h=0$ on $\Omega \backslash A$ it is easy to see that this $q$-Markov chain never escapes $A$ and so can be interpreted as a way to \emph{condition the $p$-chain to stay in $A$}.

In the case of a one-dimensional random walk, in many situations there is a unique positive harmonic function on $ \mathbb{Z}_{>0}$ (up to a multiplicative constant) which vanishes on $ \mathbb{Z}_{\leq 0}$, see \cite[Theorem 1]{Don98}. This function, properly normalized, is called the \emph{renewal function} and can be used to construct the random walk conditioned to stay positive forever, even when this event has zero probability, see \cite{BD94}. This principle will be applied below to construct $\widetilde{\Map}^{(\infty)}$ by $h$-transforming a Markovian exploration of ${\Map}^{(\infty)}$. But first we show that the process $(E_{n}: n \geq 0)$ under $ \mathbb{P}^{( \mathrm{free})}$ is itself a $h$-transformation of a killed random walk.

\subsubsection{The step distribution $\mu$.}  Recall the measure $\nu$ from \eqref{eq:nu} and let us  introduce the following measure $\mu$ on $\mathbb{Z}$ which is obtained from $\nu$, roughly speaking, by transforming half of the negative jumps into jumps of $-1$. More precisely, we put 
	$$ \left\{ \begin{array}{ll}
	\mu(2\ell) = \nu(\ell) & \mbox{ for }\ell \geq 0,\\
	\mu(-1) =  \frac{1}{2}\nu( \mathbb{Z}_{<0}),&\\ 
	\mu(-2 \ell) = \frac{1}{2}\nu(-\ell) & \mbox{ for } \ell > 0.
	\end{array} \right.$$
	Also let us introduce the function $H^\downarrow$ by setting
	\begin{equation}\label{eq:Hdowndef}
	H^{\downarrow}(\ell) \coloneqq \mathbb{P}^{( \mathrm{free})}(2|\partial  \mathfrak{m}| \geq \ell ) =  \frac{1}{ \mathrm{W}_{c}} \sum_{2j \geq \ell} W^{(j)} c_{ \mathbf{q}}^{-j}=  \frac{ \nu\left(\left(-\infty, -1-\ell/2\right]\right)}{\nu( \mathbb{Z}_{<0})}
	\end{equation} 
	for $\ell \geq 0$ and $H^\downarrow(\ell) = 0$ for $\ell < 0$, where we recall that $|\partial \map|$ is the half-perimeter of the boundary of $\map$.
	\begin{proposition} \label{prop:lawexposedboundary} If $(E_{n}: n\geq 0)$ is the length of the exposed part during the exploration of $\map$ with algorithm $ \mathcal{A}_{ \mathrm{core}}$ under $ \mathbb{P}^{( \mathrm{free})}$ then $(E_n: n \geq 0)$ is a Markov chain obtained by a $h$-transform of the random walk with i.i.d.~increments of law $\mu$ with respect to the harmonic function $H^{\downarrow}$, started from $1$ and killed on $\mathbb{Z}_{ \leq 0}$.
	\end{proposition}
	\noindent \textbf{Proof.} Recall from Section \ref{sec:freeboltz} that conditionally on the past exploration up to time $n$, if $E_{n}$ is the length of the exposed boundary then the remaining map has the law of $ \mathbb{P}^{( \mathrm{free})}$ conditioned on having perimeter at least $E_{n}$. In particular, the conditional probability to perform an event $ \mathsf{C}_{k}$ at the next peeling step is equal to 
	$$ q_k c_{ \mathbf{q}}^{k-1}\frac{\sum_{2\ell \geq E_n} W^{(\ell+k-1)} c_{ \mathbf{q}}^{-\ell-k+1}}{\sum_{2\ell \geq E_{n}} W^{(\ell)} c_{ \mathbf{q}}^{-\ell}} =  \mu(2k-2) \frac{H^{\downarrow}(E_{n}+2(k-1))}{H^{\downarrow}(E_{n})}.$$
	We recognize the form \eqref{eq:htransform} of the $H^{\downarrow}$-transform of $\mu$-random walk (killed on $ \mathbb{Z}_{\leq 0}$) in this case. Similarly, the conditional probability to perform a peeling step which identifies the peeled edge with another edge not located on the exposed boundary is equal to
	$$c_{ \mathbf{q}}^{-1}\frac{\sum_{k\geq 0} W^{(k)}c_{\mathbf{q}}^{-k} \sum_{2j \geq E_n-1} W^{(j)}c_{ \mathbf{q}}^{-j}}{\sum_{2{\ell} \geq E_{n}} W^{(\ell)} c_{ \mathbf{q}}^{-\ell}}=  \sum_{k=0}^{\infty} W^{(k)}c_{ \mathbf{q}}^{-k-1}  \frac{H^{\downarrow}(E_{n}-1)}{H^{\downarrow}(E_{n})}  =   \mu(-1)\frac{H^{\downarrow}(E_{n}-1)}{H^{\downarrow}(E_{n})} .$$
	As above, we see the $H^\downarrow$-transform of the $\mu$-walk appearing. The last case when the peeled edge is identified with another edge of the exposed boundary is similar and the conditional probability that $E_{n+1} = E_n-2k$ with $1\leq k \leq \lfloor E_n/2 \rfloor$ is, as expected, equal to 
	$$c_{\mathbf{q}}^{-1}  \frac{W^{(k-1)} c_{ \mathbf{q}}^{-k+1}\sum_{2j\geq E_n-2k} W^{(j)}c_{\mathbf{q}}^{-j} }{\sum_{2{\ell} \geq E_{n}} W^{(\ell)} c_{ \mathbf{q}}^{-\ell}} = \mu(-2k) \frac{H^\downarrow(E_n-2k)}{H^\downarrow(E_n)}.$$
	\qed \medskip

	In particular the function $H^{\downarrow}$ is bounded, equals $1$ at $0$, tends to $0$ at $\infty$, vanishes on $ \mathbb{Z}_{<0}$ and is harmonic for the $\mu$-walk on $ \mathbb{Z}_{\geq 0}$. By the optional stopping theorem we deduce that $H^{\downarrow}(k)$ is equal to the probability that a $\mu$-walk started from $k$ hits $ \mathbb{Z}_{\leq 0}$ first at $0$ in finite time. We then form the function   \begin{eqnarray}H^{\uparrow}(\ell) = H^{\downarrow}(0) + \dots + H^{\downarrow}(\ell-1)\label{eq:Hupdef}  \end{eqnarray} which is called the \emph{renewal} function of the $\mu$-walk, see e.g.~\cite{BD94}. \medskip

	Let us consider a random walk $(R_{n}:n \geq 0)$ started from $1$ with i.i.d.~increments of law $\mu$. 
	By construction $(R_{n}:n \geq 0)$ is stochastically larger than twice the $\nu$-random walk and since a $\nu$-random walk oscillates when $ \mathbf{q}$ is critical \cite[Proposition~4]{Bud15}, it follows that $R$ either oscillates or drifts towards $+\infty$. If $ \mathbf{q}$ is subcritical then a $\nu$-random walk drifts towards $-\infty$ \cite[Proposition~4]{Bud15} and it may (and typically it will) be the case that $R$ also drifts towards $-\infty$ but we excluded this case here (recall Remark \ref{rem:subcritical}). \medskip

	\noindent \textbf{Proof of Proposition \ref{prop:ratiosimple}.} As above we write   $(E_{n}: 0 \leq n \leq \tau)$ for the length of the exposed boundary when performing the peeling using algorithm $ \mathcal{A}_{ \mathrm{core}}$ under the measure $ \mathbb{P}^{( \mathrm{free})}$ and $D$ is its number of $-1$ steps until the first time $\tau \geq 1$ where $E_\tau=0$. By \eqref{eq:coreperimeter} and \eqref{eq:perimetercore}, the proposition is proved if we can show that 
	$$ \frac{\mathbb{P}^{(\mathrm{free})}(D=n+2)}{ \mathbb{P}^{(\mathrm{free})}(D=n)} \to 1, \quad \mbox{ as } n \to \infty \quad ( \mbox{along odd values}).$$ This will eventually follow from a variation on the strong ratio limit theorem \cite[Theorem 4.7]{abraham2018critical}, and in particular from the arguments of the elegant proof of Neveu \cite{Nev63}. To connect the above goal to a random walk problem, we use the fact that $(E_n)$ is the $H^\downarrow$-transform of the walk $(R_n)$ and we apply Feller's cyclic lemma (see e.g.~\cite[Chap. XII-6, Lemma 1, p412]{Fel71}).
	Using the short hand notation $\mathcal{D}_{m} = \sum_{i=0}^{m-1}  \mathbf{1}_{\Delta R_{i} = -1}$, we deduce that for $n \geq 1$ and odd we have
	 \begin{eqnarray} & &\mathbb{P}^{( \mathrm{free})}(D=n, \tau=m)\nonumber \\  &\underset{ \mathrm{h-transform}}{=}&  \mathbb{P}^{( \mathrm{free})}( \map \ne \dagger) \cdot \frac{H^\downarrow(0)}{H^\downarrow(1)} \cdot  \mathbb{P}( R_m=0, \mathcal{D}_m=n\  \& \ R_i \geq 1, \forall 0 \leq i \leq m-1)\nonumber \\
	 &=& \mathbb{P}( R_m=0, \mathcal{D}_m=n\  \& \ R_i \geq 1, \forall 0 \leq i \leq m-1)\nonumber \\
	 & \underset{ \mathrm{Cyclic\ lemma}}{=} &	  \frac{1}{m} \mathbb{P}\left(R_{m}=0,  \mathcal{D}_{m} = n\right).   \label{eq:cyclicq}\end{eqnarray}
	 Recall that $\dagger$ denotes the vertex map.
	The rough idea is then to argue that $\mathbb{P}^{( \mathrm{free})}(D=n) \approx  \mathbb{P}^{( \mathrm{free})}(D=n, \tau=N)$ with $N = [n/\mu(-1)]$ and then use the strong ratio limit theorem on the multidimensional random walk $(R_{n},  \mathcal{D}_{n})$ started at $(1,0)$ whose i.i.d.~increments are $(\Delta R_{n}, \mathbf{1}_{\Delta R_{n}=-1})$. More precisely, we will prove in Lemma \ref{lem:roughcontrol} below that for every $ \varepsilon >0$, there exists $\delta >0$, such that if $|n/m - \mu(-1)| \leq \delta$ and if $n \geq 1 /\delta$ then we have both 	  \begin{eqnarray}  \mathbb{P}( \mathcal{D}_{m}=n, R_{m}=0) &\geq&  \mathrm{e}^{- \varepsilon n}, \label{eq:gros}\\
	  \left|\frac{\mathbb{P}( \mathcal{D}_{m}=n+2, R_{m}=0)}{\mathbb{P}( \mathcal{D}_{m}=n, R_{m}=0)} -1 \right| &\leq& \varepsilon, \label{eq:ratio}	  \end{eqnarray}
	where implicitly in the above display and in the rest of the proof, we restrict to values of $n,m$ for which the probabilities considered are positive. Indeed, the random walk $(R, \mathcal{D})$ is not aperiodic since $R_n+ \mathcal{D}_n$ is always odd, hence the $n+2$. To come back to $ \mathbb{P}(D=n)$, notice that since $ \mathcal{D}_{m}$ is distributed as $ \mathrm{Binomial}(m, \mu(-1))$, easy large deviation estimates show that for the $\delta$ fixed above, for all $m$ large enough we have,
	$$ \mathbb{P}\left(\left| \frac{ \mathcal{D}_{m}}{m} - \mu(-1)\right| \geq \delta \right) \leq  \mathrm{e}^{- c m},$$ for some $c>0$. Summing over all $ m \geq n$ such that   $|n/m- \mu(-1)| \geq \delta$, we deduce using \eqref{eq:cyclicq} that  
	 \begin{eqnarray} \label{eq:sumpresque} \Big|\mathbb{P}^{( \mathrm{free})}(D=n) -  \underbrace{\sum_{\begin{subarray}{c} m \geq n \\ |\frac{n}{m}- \mu(-1)| \leq \delta \end{subarray}} \frac{1}{m}\mathbb{P}\left( \mathcal{D}_m=n \mbox{ and } R_m=0 \right)}_{ \displaystyle =: \ \Sigma(n)} \Big| \leq \sum_{m \geq n} \mathrm{e}^{-c m} \leq  c' \cdot \mathrm{e}^{-c n}, \end{eqnarray} for some $c' >0$. In words, the probability $\mathbb{P}^{( \mathrm{free})}(D=n)$ is approximated up to an exponentially small probability by the sum $\Sigma (n)$ in the last display.  Thanks to \eqref{eq:ratio}, $|\Sigma(n+2)/\Sigma(n)-1| \leq \varepsilon$ since the ratio of each summand is close to $1$ (we are neglecting boundary effects here). On the other hand, using \eqref{eq:cyclicq} and \eqref{eq:gros} when $ \varepsilon \to 0$, we deduce that $ (\mathbb{P}^{( \mathrm{free})} (D=n))^{1/n} \to 1$ as $n \to \infty$. Combining this with \eqref{eq:sumpresque} we deduce that $\mathbb{P}^{( \mathrm{free})} (D=n) \sim \Sigma(n)$ as $n \to \infty$. Consequently, the ratio $\mathbb{P}^{( \mathrm{free})} (D=n+2)/\mathbb{P}^{( \mathrm{free})} (D=n)$ also belongs to $(1- 2 \varepsilon,  1+ 2 \varepsilon)$ asymptotically.\qed \medskip 
	
	It remains to prove the lemma we used in the course of the proof. 
	\begin{lemma}  \label{lem:roughcontrol}  	For every $ \varepsilon >0$, there exists $\delta >0$, such that if $|n/m - \mu(-1)| \leq \delta$ and if $n \geq 1 /\delta$ then we have both 
	  \begin{eqnarray}  \mathbb{P}( \mathcal{D}_{m}=n, R_{m}=0) &\geq&  \mathrm{e}^{- \varepsilon n}, \\
	  \left|\frac{\mathbb{P}( \mathcal{D}_{m}=n+2, R_{m}=0)}{\mathbb{P}( \mathcal{D}_{m}=n, R_{m}=0)} -1 \right| &\leq& \varepsilon,
	  \end{eqnarray}
	  where we restrict to values of $n,m$ for which the probabilities are positive.
	\end{lemma}
	\noindent \textbf{Proof:} Fix $ \varepsilon>0$. The first point consists in showing that $\mathbb{P}( \mathcal{D}_{m}=n, R_{m}=0)$ does not decay exponentially fast. For this, it suffices to exhibit a scenario where  $\mathcal{D}_{m}=n$ and  $R_{m}=0$ which happens with a not too small probability. Notice first that since $ \mathcal{D}_m$ has a $ \mathrm{Binomial}( m,\mu(-1))$ distribution, we can find $\delta >0$ small enough so that if $m,n$ satisfy the requirements of the lemma then 
	$$ \mathbb{P}( \mathcal{D}_{m-1} =n) \geq \mathrm{e}^{ - \varepsilon n}. $$ By restricting furthermore to paths having a bound on their $R$-increments, we can find $A >0$ large enough so that up to decreasing $\delta>0$ we have 
	$$ \mathbb{P}( D_{m-1} = n \mbox{ and } |R_{i+1}-R_i| \leq A, \forall 0 \leq m-2) \geq  \mathrm{e}^{-  \varepsilon n},$$
	as soon as $|n/m- \mu(-1)| \leq \delta$ and $n \geq 1/\delta$.
	Now, from \cite[Eq (5.13)]{CurStFlour} the measure $\nu$, hence $\mu$, always has a polynomial tail on the left, namely~$\mu(-k) \geq \mathrm{c} k^{-5/2}$ asymptotically as $k \to \infty$ for some $c>0$. Hence, in the above scenario, since $R_{m-1} \leq Am$, a single large negative jump of $R$ at time $m$ could yield to the value $(n,0)$ with the additional cost of  a polynomially decaying probability. This proves the first point of the lemma.\\
	The second point of the lemma follows from the first point combined with Neveu's proof of the strong ratio limit theorem which we now recall. Fix a possible increment $(a,b) \in \mathbb{Z}^2$ of the walk $( \mathcal{D},R)$. Following Neveu \cite{Nev63} we write $\mathbb{P}( (\Delta  \mathcal{D}_0, \Delta R_0)= (a,b) \mid ( \mathcal{D}_m,R_m)=(n,0))$ in two different ways. Introducing $N_m$ the number of increments of $( \mathcal{D},R)$ equal to $(a,b)$ up to time $m$ and using the permutation symmetry of the increments:
 \begin{eqnarray} \mathbb{E}\left.\left[ \frac{N_{m}}{m} \,\right| ( \mathcal{D}_m,R_m)=(n,0)\right] &\underset{ \mathrm{symmetry}}{=}& \mathbb{P}( (\Delta  \mathcal{D}_0, \Delta R_0)= (a,b) \mid ( \mathcal{D}_m,R_m)=(n,0)) \label{eq:neveu}\\ & \underset{ \mathrm{Markov}}{=} & \mathbb{P}((\Delta  \mathcal{D}_0, \Delta R_0)= (a,b)) \cdot \frac{ \mathbb{P}(( \mathcal{D}_{m-1},R_{m-1})=(n-a,-b))}{ \mathbb{P}(( \mathcal{D}_m,R_m)=(n,0))}.\nonumber    \end{eqnarray}
 Given the periodicity conditions of the walk $( \mathcal{D},R)$, to prove the lemma it suffices to show that for any $ \eta>0$ and for any $(a,b)$ in the support of the increment, there exists $\delta >0$ such that the ratio in the last display lies in $[1- \eta, 1+ \eta]$ as soon as $|n/m - \mu(-1)| \leq \delta$ and $n \geq 1 /\delta$. To see this, we examine the left-hand side of \eqref{eq:neveu}. Indeed, since $N_m$ has distribution $ \mathrm{Binomial}( m, p_{a,b})$ with $p_{a,b}=\mathbb{P}((\Delta  \mathcal{D}_0, \Delta R_0)= (a,b))$, an easy large deviation estimate shows that for all $m$ large enough
 $$ \mathbb{P}\left(\left| \frac{N_{m}}{ p_{a,b}\cdot m} -  1\right| \geq \eta\right) \leq \mathrm{e}^{- c_\eta m },$$ for some constant $c_\eta>0$. Applying this estimate to \eqref{eq:neveu} we find that 
 \begin{align*}
 \left| \frac{ \mathbb{P}(( \mathcal{D}_{m-1},R_{m-1})=(n-a,-b))}{ \mathbb{P}(( \mathcal{D}_m,R_m)=(n,0))}-1 \right| &= \left| \mathbb{E}\left.\left[ \frac{N_{m}}{ p_{a,b} \cdot m}  \,\right| ( \mathcal{D}_m,R_m)=(n,0)\right]-1 \right| \\
 &\leq \eta +  \frac{1}{p_{a,b}}\cdot \frac{ \mathrm{e}^{-c_\eta m }}{\mathbb{P}(( \mathcal{D}_m,R_m)=(n,0))}.\end{align*}
 Using the first point of the lemma, for fixed $\eta>0$, we can decrease $\delta$ so that $\mathbb{P}(( \mathcal{D}_m,R_m)=(n,0))$ is asymptotically much larger than $\mathrm{e}^{-c_\eta m }$, thus asymptotically bounding the above display by $2 \eta$. This finishes the proof of the lemma.\qed \medskip

	 \section{Infinite Simple Boltzmann maps of the half-plane}
	 In this section, we introduce the infinite simple Boltzmann maps of the half-plane $ \hat{\Map}^{(\infty)}$ of law $ \hat{ \mathbb{P}}^{(\infty)}$ which will be the local limit of the Boltzmann maps with a large simple boundary (Theorem \ref{thm:localsimple}). To do so, we first construct $\widetilde{\Map}^{(\infty)}$, the infinite Boltzmann map of the half-plane (with a general boundary) conditioned on having an infinite core. Then we simply define $ \hat{\Map}^{(\infty)}$ as the infinite core of $\widetilde{\Map}^{(\infty)}$ (Definition \ref{def:Minftyhat}).
	 
	  In the case when $\gulp<\infty$ (the so-called dilute case), $\hat\Map^{(\infty)}$ can  be directly defined as the unique infinite simple component of $\Map^{(\infty)}$. Equivalently, $\widetilde{\Map}^{(\infty)}$ is obtained by conditioning $\Map^{(\infty)}$ on the event of positive probability that the root edge is on the infinite simple component. However, when $\gulp=\infty$, the last event has probability zero and we shall use a $h$-transformation to construct $\widetilde{\Map}^{(\infty)}$.

	 \subsection{Construction of $\widetilde{\Map}^{(\infty)}$}
	 
\label{sec:Mtilde}

Compared to  early works on random planar maps \cite{AS03}, the infinite random map $\widetilde{\Map}^{(\infty)}$ will not be defined as a local limit of finite random maps. It is rather constructed directly in the infinite setting as in \cite{AR13,CurPSHIT}, and this will be done in three steps:
\begin{enumerate}
\item \emph{Uniqueness.} We shall state in Proposition \ref{prop:characSHP} a spatial Markov property that uniquely characterizes the law of $\widetilde{\Map}^{(\infty)}$ provided it exists.
\item \emph{Peeling process.} Assuming its existence, we prove that a random map $\widetilde{\Map}^{(\infty)}$ satisfying the above spatial Markov property obeys a peeling process whose transitions are explicit and indeed sum to $1$ (Lemma \ref{lem:sum1}).
\item \emph{Construction.} We finally construct $\widetilde{\Map}^{(\infty)}$ as the increasing union of maps whose growing mechanism is given by the peeling process driven by a well-chosen algorithm (Proposition \ref{prop:constructSUIHPQ}). We check that it satisfies the desired spatial Markov property.
\end{enumerate}
See \cite{CurStFlour} for several applications of this strategy to construct infinite random maps.

\begin{figure}[!h]
	 \begin{center}
	 \includegraphics[width=7cm]{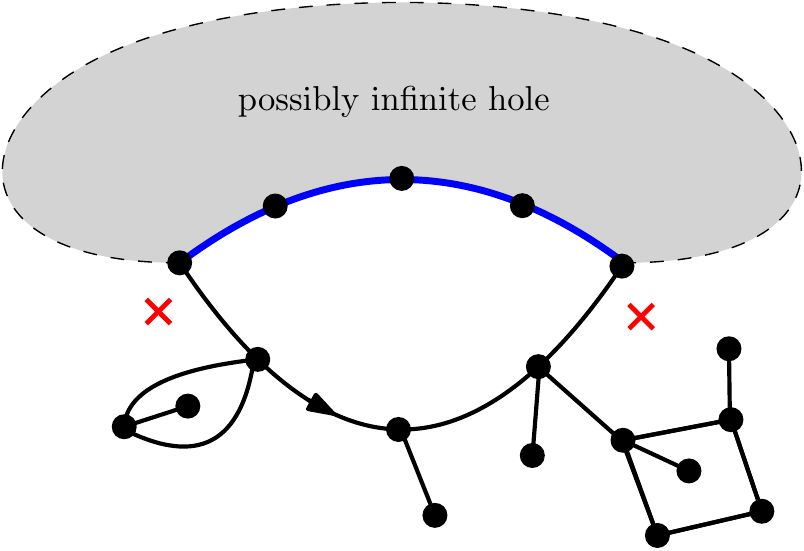}
	 \caption{ \label{fig:goodsubmap}A good submap $ \tilde{ \mathfrak{e}}$ (which may be a submap  of the half-plane). The exposed boundary is indicated by the thick blue line. Notice that no dangling parts are attached to the boundary points of the exposed part (marked by red crosses). Here $| \partial^{*} \tilde{ \mathfrak{e}}| - | \partial \tilde{ \mathfrak{e}}| = 2 - 9=-7$. }
	 \end{center}
	 \end{figure}

\paragraph{Good submaps and explorations.} A map $\map$ of the half-plane is an infinite planar map with one end such that the root face has an infinite degree. Recall from Section \ref{sec:freeboltz} the notion of submap $ \tilde{ \mathfrak{e}}$ where the perimeter of the unique hole is \emph{a priori} not determined.  Clearly the notion $ \tilde{ \mathfrak{e}} \Subset  \mathfrak{m}$ makes sense for maps $ \mathfrak{m}$ of the half-plane. In what follows we shall restrict to those submaps $\tilde{ \mathfrak{e}}$ which can be obtained by a filled-in peeling algorithm that always peel on the exposed boundary\footnote{We do so because we will use an $h$-transformation on the length of the exposed boundary. Moreover, in a submap with undetermined perimeter, the algorithm "does not know" the perimeter of the dotted line of Figure \ref{fig:submaprandom}, so to stay on solid ground it is better to stick to algorithms that always peel on the exposed boundary. }. Those \emph{good} submaps are characterized as follows: $\tilde{\mathfrak{e}}$ is a \emph{good} submap if the simple component carrying the root edge shares at least one edge with the exposed boundary and if all other simple components of $\tilde{\mathfrak{e}}$ do not share a vertex with the exposed boundary (see Figure \ref{fig:goodsubmap}). For example, the  submaps explored using algorithm $\mathcal{A}_{ \mathrm{core}}$ (until its stopping time) are all of this type. 
	
We take the quantity $|\partial^* \tilde{ \mathfrak{e}}| - |\partial \tilde{ \mathfrak{e}}|$ to be the difference between the \emph{half}-perimeters of the exposed and internal boundaries. Recall the definition of $H^\uparrow$ as introduced in \eqref{eq:Hupdef}. 

	\begin{proposition} \label{prop:characSHP} If it exists, there is a unique law $ \widetilde{ \mathbb{P}}^{(\infty)}$ supported by infinite maps of the half-plane (with general boundary) so that its core (the simple component containing the root edge) has almost surely infinite perimeter and whose law is characterized by 
	 \begin{eqnarray} \label{eq:characSHP} \widetilde{\mathbb{P}}^{(\infty)}(   \tilde{\mathfrak{e}} \Subset \map) = w_{ \mathbf{q}}(  \tilde{\mathfrak{e}}) \cdot  {c_{ \mathbf{q}}}^{ | \partial^{*} \tilde{ \mathfrak{e}}| - | \partial \tilde{ \mathfrak{e}}|} \cdot H^\uparrow( p_{\mathrm{exposed}})  \end{eqnarray}
	for any good submap $  \tilde{\mathfrak{e}}$ with  a unique hole and such that the perimeter of the exposed boundary is $ p_{ \mathrm{exposed}}$. 
	\end{proposition}
	\noindent \textbf{Proof.} Let $\widetilde{\Map}^{(\infty)}$ be a random map of the half-plane of law $ \widetilde{\mathbb{P}}^{(\infty)}$. The law of $\widetilde{\Map}^{(\infty)}$ is fully characterized by the law of balls of radius $r$ around the root edge for $r \geq 1$. But since by hypothesis $\widetilde{\Map}^{(\infty)}$ is one-ended almost surely, the last laws are completely characterized by the probabilities $\widetilde{\mathbb{P}}^{(\infty)}(  \tilde{\mathfrak{e}} \Subset \map)$ for finite good submaps $ \tilde{ \mathfrak{e}}$ with a unique hole as involved in the proposition.
	 \qed \medskip 
	
	Following the strategy sketched at the beginning of this section, if $\widetilde{\Map}^{(\infty)}$ is a map of law $\widetilde{ \mathbb{P}}^{(\infty)}$ (provided it exists), it is easy to deduce the form of the probability transition for any \emph{good} peeling algorithm. 
	
	\begin{proposition}\label{prop:peelinghalfplaneMC} Let us assume that a law $\widetilde{ \mathbb{P}}^{(\infty)}$ satisfying \eqref{eq:characSHP} does exist. Fix a good submap $ \tilde{ \mathfrak{e}}$ and an edge $a$ of its exposed boundary. Under $\widetilde{ \mathbb{P}}^{(\infty)}$ conditionally on $\{ \tilde{ \mathfrak{e}} \Subset \map\}$, the peeling of the edge $a$ yields the events $\mathsf{C}_{k}$, $\mathsf{G}_{k,\infty}$ and $\mathsf{G}_{\infty, k}$ for $k \geq 0$ with respective probabilities
	\begin{equation}\label{eq:Mtildetransitions}
	q_{k} c_{ \qseq}^{k-1} \frac{H^{\uparrow}(p_{ \mathrm{new}})}{H^{\uparrow}(p_{ \mathrm{old}})}, \quad W^{(k)} c_{ \qseq}^{-k-1} \frac{H^{\uparrow}(p_{ \mathrm{new}})}{H^{\uparrow}(p_{ \mathrm{old}})}, \quad \mbox{and} \quad W^{(k)} c_{ \qseq}^{-k-1} \frac{H^{\uparrow}(p_{ \mathrm{new}})}{H^{\uparrow}(p_{ \mathrm{old}})}
	\end{equation} where $p_{ \mathrm{old}}$ and $p_{ \mathrm{new}}$ are respectively the length of the exposed boundary before and after the peeling step. Conditionally on these events the possible finite hole created is filled-in with an independent $ \mathbf{q}$-Boltzmann map with the correct perimeter.
\end{proposition}
\noindent \textbf{Proof.} The event $  \mathsf{C}_k$ happens if and only if $ \tilde{ \mathfrak{e}} ' \Subset \map$ where $\tilde{ \mathfrak{e}} '$ is the map obtained from $ \tilde{ \mathfrak{e}}$ by gluing a face of degree $2k$ on $a$. The internal boundaries of $ \tilde{ \mathfrak{e}}$ and $ \tilde{ \mathfrak{e}}'$ coincide and  we have $|\partial^* \tilde{ \mathfrak{e}}'| =|\partial^* \tilde{ \mathfrak{e}}'|+k-1$. By \eqref{eq:characSHP} the conditional probability of $ \mathsf{C}_k$ is then 
$$ \frac{\widetilde{\mathbb{P}}^{(\infty)}(   \tilde{\mathfrak{e}}' \Subset \map)}{\widetilde{\mathbb{P}}^{(\infty)}(   \tilde{\mathfrak{e}} \Subset \map)} = q_k \cdot c_{ \mathbf{q}}^{k-1} \frac{H^\uparrow(2 |\partial^* \tilde{ \mathfrak{e}}'|)}{H^\uparrow(2 |\partial^* \tilde{ \mathfrak{e}}|)}.$$
The other transitions are computed similarly. \qed \medskip
	
	 We need now to check that the numbers appearing in the last proposition indeed define probability transitions (unconditionally on the existence of $\widetilde{ \mathbb{P}}^{(\infty)}$):
	\begin{lemma} \label{lem:sum1}The  probability transitions of \eqref{eq:Mtildetransitions} sum up to $1$.
	\end{lemma}
	\noindent \textbf{Proof.} Let us consider the situation in which the exposed boundary has length $p\geq 1$ and the peeling algorithm selects the $\ell$th edge from the left, with $1\leq \ell \leq p$. 
	Expressing the probabilities above in terms of the law $\nu$ we need to check that
	\begin{equation}\label{eq:genharmonicHup}
	\sum_{k\geq 0} \nu(k) H^\uparrow(p+2k) + \sum_{k \leq -1} \frac{1}{2}\nu(k)H^\uparrow((p+2k) \vee (p-\ell)) + \sum_{k \leq -1} \frac{1}{2}\nu(k)H^\uparrow((p+2k) \vee (\ell-1)) = H^\uparrow(p).
	\end{equation}
	We have already seen in Proposition \ref{prop:lawexposedboundary} that $H^{\downarrow}$ is the pre-renewal function for the $\mu$-walk and is thus harmonic. Since the $\mu$-walk is not drifting towards $-\infty$ it follows that its renewal function $ H^{\uparrow}$ is also harmonic for the $\mu$-walk, see \cite[Appendix B]{CC08}. In our context, this means that the above probability transitions indeed sum-up to $1$ \emph{in the particular case of the algorithm $ \mathcal{A}_{ \mathrm{core}}$} which peels the left-most edge of the exposed boundary. 
	Therefore \eqref{eq:genharmonicHup} is satisfied for $\ell=1$.
	
	Next we compare the left-hand side of \eqref{eq:genharmonicHup} when peeling the $(\ell+1)$th edge to the case where the $\ell$th edge is peeled. 
	If $1\leq \ell \leq p-1$ the difference in the second term is
	\begin{align*}
	\sum_{k \leq -1} &\frac{1}{2}\nu(k)(H^\uparrow((p+2k) \vee (p-\ell-1)) - H^\uparrow((p+2k) \vee (p-\ell))) \\
	&= \sum_{2k \leq -\ell-1} \frac{1}{2}\nu(k)(H^\uparrow(p-\ell-1) - H^\uparrow(p-\ell)) \\
	&\stackrel{\eqref{eq:Hupdef}}{=}-\sum_{2k \leq -\ell-1} \frac{1}{2}\nu(k) H^\downarrow(p-\ell-1) \stackrel{\eqref{eq:Hdowndef}}{=} -\frac{1}{2}\nu( \mathbb{Z}_{<0}) H^\downarrow(\ell-1)H^\downarrow(p-\ell-1).
	\end{align*}
	It is canceled by the difference in the third term, which evaluates to $\frac{1}{2}\nu( \mathbb{Z}_{<0}) H^\downarrow(p-\ell-1)H^\downarrow(\ell-1)$. By induction the identity \eqref{eq:genharmonicHup} follows for all $\ell$.\qed \medskip

Given the last result, for any filled-in good lazy-peeling algorithm $ \mathcal{A}$, the transitions \eqref{eq:Mtildetransitions} for the events $ \mathsf{C}_{\cdot}$ or $ \mathsf{G}_{\cdot, \cdot}$ \emph{define} for us a Markov chain $ (\tilde{\mathbf{E}}_n : n \geq 0)$ of growing good submaps starting from the initial edge map and which never stops (since the length of the exposed boundary never drops to $0$). We shall construct the law $ \widetilde{ \mathbb{P}}^{(\infty)}$ using this chain for a particular algorithm which ensures that we indeed build a map of the half-plane (as opposed to an algorithm that peels only a subregion of the space without discovering the neighborhood of the origin).

\begin{center}
	\fbox{\begin{minipage}{15cm}
			\paragraph{Algorithm $\mathcal{A}_{\text{metric}}$:}  For any good submap $ \tilde{ \mathfrak{e}}$  with a unique hole, let $ \mathcal{A}_{ \mathrm{metric}}( \tilde{ \mathfrak{e}})$ be the left-most edge of the exposed boundary with an endpoint minimizing the graph  distance, inside $ \tilde{ \mathfrak{e}}$, to the origin of the root edge. We then fill-in the holes of finite perimeter we may create on the way.
		\end{minipage}}
	\end{center}
	
	\begin{proposition} \label{prop:constructSUIHPQ}  Let $ (\tilde{\mathbf{E}}_n : n \geq 0)$ be the Markov chain of growing good submaps whose probability transitions are given by Proposition \ref{prop:peelinghalfplaneMC} with the peeling algorithm $ \mathcal{A}_{ \mathrm{metric}}$. Then the map
	$$  \widetilde{\Map}^{(\infty)}:=\bigcup_{n \geq 0} \tilde{ \mathbf{E}}_n,$$
	is a random infinite map of the half-plane  
	which satisfies \eqref{eq:characSHP}.
	\end{proposition}
	\noindent \textbf{Proof.} The proof is similar to  \cite[Section 1.3]{CurPSHIT} or \cite[Proposition 6.5]{CurStFlour}. There are two non-trivial points in the proposition. First, one needs to prove that $\widetilde{\Map}^{(\infty)}$ as defined above is indeed a map of the half-plane and second that it satisfies \eqref{eq:characSHP}. For the first point, the problem that could appear is that some vertex $x$ remains exposed on $ \partial^{*} \tilde{\mathbf{E}}_{n}$ forever (i.e.\ is never swallowed by the process). This cannot happen a.s., since if the perimeter of the exposed boundary is $p \geq 3$, the next step the Markov chain may swallow the point on the right or on the left of the peeled edge  via an event $ \mathsf{G}_{0, \infty}$ or $ \mathsf{G}_{\infty,0}$ (i.e.~this vertex becomes an internal vertex of $ \tilde{ \mathbf{E}}_{n+1}$). During such an event the length of the exposed boundary may decrease by one or two so this happens with probability at least
	$$ \inf_{p \geq 3 }\frac{H^{\uparrow}(p-1) \wedge H^{\uparrow}(p-2)}{H^{\uparrow}(p)} c_{ \mathbf{q}}^{-1}>0.$$  If $p=1$ or $2$ then such a vertex can be swallowed in two steps with some probability $c>0$. It easily follows from the definition of $ \mathcal{A}_{ \mathrm{metric}}$ that the minimal distance to the origin of a point of $\partial^*    \tilde{\mathbf{E}}_n$ tends to $\infty$ almost surely and so $ \widetilde{\Map}^{(\infty)}$ is indeed a map of the half-plane almost surely.\\
	We then need to check that $\widetilde{\Map}^{(\infty)}$ satisfies \eqref{eq:characSHP}. Fix a good submap $ \tilde{e}_{n_{0}}$ of the half-plane with a unique hole of infinite perimeter which can be obtained as the result of $n_{0}$ steps of a good filled-in lazy-exploration $\tilde{e}_0 \Subset \cdots \Subset \tilde{e}_{n_0}$.  Our goal is then to prove that 
	  \begin{eqnarray} \label{eq:goalconstruction} \widetilde{\mathbb{P}}^{(\infty)}( \tilde{e}_{n_{0}} \Subset\widetilde{\Map}^{(\infty)} ) = w_{ \mathbf{q}}( \tilde{e}_{n_{0}}) c_{ \mathbf{q}}^{ |\partial^{*} \tilde{e}_{n_{0}}|- | \partial \tilde{ e}_{n_{0}}|} H^\uparrow(p_{\text{exposed}}),  \end{eqnarray}
	  where $p_{\text{exposed}}$ is the length of the exposed boundary of $\tilde{e}_{n_0}$.
	  To do this, we first design another good filled-in peeling algorithm $ \mathcal{A}'_{\mathrm{metric}}$ so that if $ \tilde{ \mathfrak{e}} = \tilde{e}_{i}$ for some  $i < n_0$ then $ \mathcal{A}_{ \mathrm{metric}}'( \tilde{ \mathfrak{e}})$ is the edge that is to be peeled to pass from $ \tilde{e}_i$ to $ \tilde{e}_{i+1}$, and otherwise  we put $\mathcal{A}_{ \mathrm{metric}}'( \tilde{ \mathfrak{e}}) = \mathcal{A}_{ \mathrm{metric}}( \tilde{ \mathfrak{e}}) $. Roughly speaking, the filled-in exploration using this algorithm first decides whether or not we have $  \tilde{ e}_{n_{0}} \Subset \map$ and then performs the  metric exploration with $ \mathcal{A}_{ \mathrm{metric}}$. We denote by $(\tilde{ \mathbf{E}}_n : n \geq 0)$ and $ (\tilde{ \mathbf{E}}'_n: n \geq 0)$  the Markov chains on growing maps ruled by the transitions \eqref{eq:Mtildetransitions} and peeling algorithm $ \mathcal{A}_{ \mathrm{metric}}$ and $ \mathcal{A}_{ \mathrm{metric}}'$ respectively. Adapting the above argument we get that $$ \widetilde{\Map}^{(\infty)\prime} = \bigcup_{n \geq 0} \tilde{\mathbf{E}}_{n}',$$ is almost surely a map of the half-plane. Using the explicit transition probabilities \eqref{eq:Mtildetransitions} we have
	$$ \mathbb{P}( \tilde{e}_{n_{0}} \Subset \widetilde{\Map}^{(\infty)\prime}) = \mathbb{P}( \tilde{ \mathbf{E}}'_{n_{0}} = \tilde{e}_{ n_{0}}) = w_{ \mathbf{q}}( \tilde{e}_{n_{0}}) c_{ \mathbf{q}}^{ |\partial^{*} \tilde{e}_{n_{0}}|- | \partial \tilde{ e}_{n_{0}}|} H^\uparrow(p_{\text{exposed}}).$$
	It remains to prove that $\widetilde{\Map}^{(\infty)} = \widetilde{\Map}^{(\infty)\prime}$ in law to deduce our goal \eqref{eq:goalconstruction}. To see this, let us introduce the stopping time $\tau$ (resp.~ $\tau'$) to be the first time at which the minimal graph distance of vertices of $\partial^*  \tilde{  \mathbf{E}}_{n}$ (resp.~$ \partial^* \tilde{ \mathbf{E}}_{n}'$) to the origin of the root edge is larger than $d$, were $d\geq 0$ is chosen to be larger than the diameter of $ \tilde{e}_{n_{0}}$. By the above consideration we have $\tau,\tau'< \infty$ almost surely and an easy extension of the last display  to the almost sure finite  stopping times $\tau,\tau'$ shows that for any $ \tilde{ e }$ with exposed boundary length $p_{\text{exposed}}$ we have 
	$$ \mathbb{P}( \tilde{ \mathbf{E}}_{\tau} = \tilde{ e}) = w_{ \mathbf{q}}( \tilde{e}) c_{ \mathbf{q}}^{ |\partial^{*} \tilde{e}|- | \partial \tilde{ e}|} H^\uparrow(p_{\text{exposed}})\mathbf{1}_{  \mbox{ $\tilde{e}$ can be obtained as $ \tilde{ \mathbf{E}}_{\tau}$}},$$
	$$\mathbb{P}( \tilde{ \mathbf{E}}'_{\tau'} = \tilde{ e}) = w_{ \mathbf{q}}( \tilde{e}) c_{ \mathbf{q}}^{ |\partial^{*} \tilde{e}|- | \partial \tilde{ e}|} H^\uparrow(p_{\text{exposed}})\mathbf{1}_{  \mbox{ $\tilde{e}$ can be obtained as $ \tilde{ \mathbf{E}}'_{\tau'}$}}.$$
	But by our choice of $d$ in the definition of $\tau$ and the properties of the algorithm $ \mathcal{A}_{ \mathrm{metric}}$ and $ \mathcal{A}_{ \mathrm{metric}}'$, it is easy to see that the possible outcomes of $ \tilde{ \mathbf{E}}_{\tau}$ or $ \tilde{ \mathbf{E}}'_{\tau'}$ are the same. This implies that $ \tilde{ \mathbf{E}}_{\tau} = \tilde{\mathbf{E}}'_{\tau'}$ in distribution and hence that we have  $ \bigcup_{n \geq 0} \tilde{ \mathbf{E}}'_{n}= \widetilde{\Map}^{(\infty)\prime} \stackrel{(d)}{=} \widetilde{\Map}^{(\infty)} = \bigcup_{n \geq 0} \tilde{ \mathbf{E}}_{n}$ as desired.
	\qed \medskip

	Interpreted in terms of the $h$-transform, we see from \eqref{eq:Mtildetransitions} that the exploration of  $\widetilde{\Map}^{(\infty)}$ with the good algorithm $ \mathcal{A}_{ \mathrm{metric}}$ is just the $H^{\uparrow}$-transform of the exploration of  ${\Map}^{(\infty)}$ with algorithm $ \mathcal{A}_{ \mathrm{metric}}$. In particular, when  $\gulp<\infty$, since a $\nu$-random walk cannot drift to $\infty$ \cite[Proposition~4]{Bud15} it must have a finite mean and so $\mu$ has a finite mean as well. Actually, in this case $\mu$ drifts to $\infty$ because $\nu$ oscillates. It follows that $H^{\uparrow}$ is bounded and converges. Furthermore one may check that 
	$$ \frac{H^{\uparrow}(1)}{H^{\uparrow}(\infty)}  \underset{\eqref{eq:Hupdef}, \eqref{eq:Hdowndef}}{=} \frac{1}{ \mathbb{E}^{( \mathrm{free})}[2 | \partial \map|]+1} = \frac{ \mathrm{W}_{c}}{2 c_{ \mathbf{q}}^{-1} \mathrm{W}'_{ \mathrm{c}} +  \mathrm{W}_{c}},$$ and this is also the probability that the exploration with $ \mathcal{A}_{ \mathrm{metric}}$ does not stop. Using \eqref{eq:characSHP} one can thus see that $ \widetilde{ \Map}^{(\infty)}$ is nothing else but $ \Map^{(\infty)}$ conditioned on the event of positive probability $ \{| \partial \mathrm{Core}( \Map^{(\infty)})|=\infty\}$.
	
	Alternatively, by invariance under translation (see \cite[Proposition 6.6]{CurStFlour}) one can define $\widetilde{\Map}^{(\infty)}$ by shifting the root edge of $\Map^{(\infty)}$ to the first edge to its right on the boundary belonging to an infinite simple component, and biasing by $1 / (1+ 2 | \partial C_{0}|)$ where $C_{0}$ is the component dangling on the left of the root edge in $\widetilde{\Map}^{(\infty)}$. 
When $\gulp=\infty$, the good explorations of $\Map^{(\infty)}$ will almost surely terminate. In this case, the strength of the $h$-transformation still enables us to interpret $\widetilde{\Map}^{(\infty)}$ as conditioning $\Map^{(\infty)}$ on the complementary zero-probability event.

\subsection{Defining $\hat{\Map}^{(\infty)}$}

We can now introduce the main character of this paper:

\begin{definition} \label{def:Minftyhat} Let $ \mathbf{q}$ be a critical weight sequence. We define $ \hat{\Map}^{(\infty)}$, the infinite Boltzmann map of the half-plane with a simple boundary, whose law is denoted by $ \hat{ \mathbb{P}}^{(\infty)}$, as the infinite core of $ \widetilde{\Map}^{(\infty)}$, i.e.
$$ \hat{\Map}^{(\infty)} \coloneqq  \mathrm{Core}(\widetilde{\Map}^{(\infty)}).$$
\end{definition}

\section{Simple peeling exploration of $\hat{\Map}^{(\infty)}$}

\begin{figure}
	  \begin{center}
	 	  \includegraphics[width=12cm]{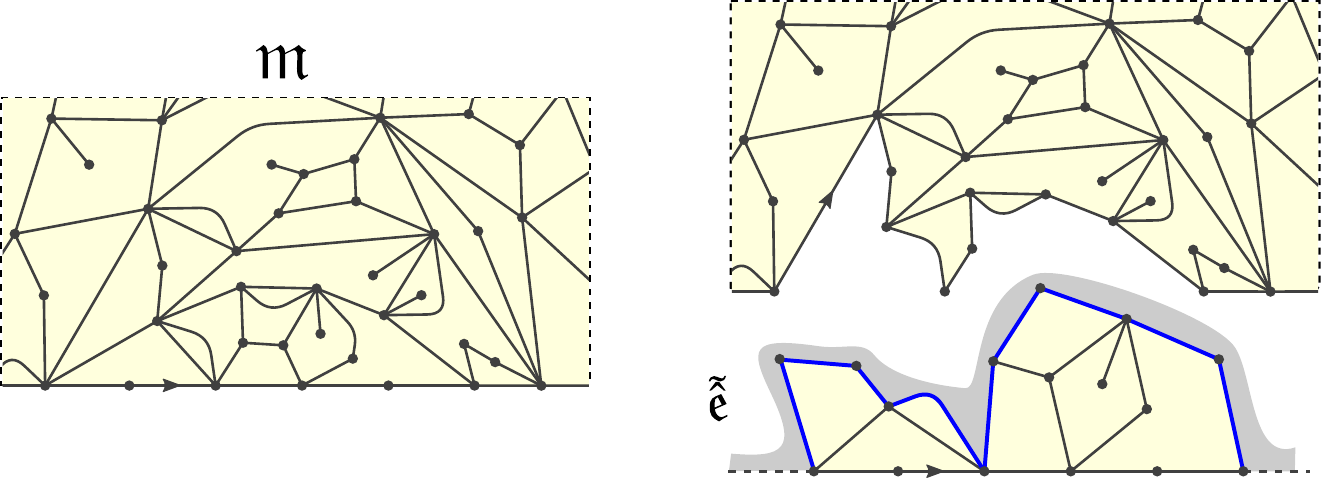}
	 	  \caption{\label{fig:simplemarkov} An example of a simple submap $\hattilde{\mathfrak{e}}$ with $|\partial^*\hattilde{\mathfrak{e}}|-|\partial\hattilde{\mathfrak{e}}|=(9-5)/2=2$. Here we have $\hattilde{ \mathfrak{e}} \hat{\Subset}  \mathfrak{m}$ where $\mathfrak{m}$ is a $\partial$-simple infinite half-planar map. The remaining map is shown on the top right.}
	 	  \end{center}
	 \end{figure}

We now turn to the study of  the filled-in simple peeling process of $ \hat{ \mathbb{P}}^{(\infty)}$. The computation of the simple peeling transitions will be obtained via the lazy peeling process and the core decomposition of $\widetilde{\Map}^{(\infty)}$. For our purpose, all we need is the following simple version of the spatial Markov property. Let $\hattilde{\mathfrak{e}}$ be a submap  with both internal and exposed boundaries simple (hence the tilde and the hat in the notation). As in the non-simple case, we write $\hattilde{ \mathfrak{e}} \hat{\Subset}  \mathfrak{m}$ if the $\partial$-simple map $ \mathfrak{m}$ can be obtained by filling-in the hole of $\hattilde{ \mathfrak{e}}$ with a $\partial$-simple map. See Figure \ref{fig:simplemarkov} for an example.
	 
	 \begin{proposition}[Spatial Markov Property under $ \hat{\mathbb{P}}^{(\infty)}$] \label{prop:spatialMarkovsimple} Let $\hattilde{\mathfrak{e}}$ be a simple submap with a unique  hole. Denote by $ | \partial^{*} \hattilde{\mathfrak{e}}|-| \partial \hattilde{\mathfrak{e}}|$ the difference of the exposed and internal boundary half-perimeters of $\hattilde{\mathfrak{e}}$. Then we have 
	 $$ \hat{\mathbb{P}}^{(\infty)}\left(\hattilde{\mathfrak{e}} \hatSubset  \map \right) =  {\hat{c}_{ \mathbf{q}}}^{| \partial^{*}\hattilde{\mathfrak{e}}|-| \partial \hattilde{\mathfrak{e}}|} \cdot w_{ \mathbf{q}}( \hattilde{\mathfrak{e}}),$$
	 with $\hat{c}_{\mathbf{q}}$ given by \eqref{def:cqhat}.
	 Furthermore, conditionally on $\hattilde{\mathfrak{e}} \hatSubset \map$ the  map filling-in the hole of  $\hattilde{\mathfrak{e}}$ has law $\hat{ \mathbb{P}}^{(\infty)}$.
	\end{proposition}
	To prove the above proposition we will rely on lazy-exploration of $\widetilde{\Map}^{(\infty)}$. We first compute the law of the core decomposition under $ \widetilde{\mathbb{P}}^{(\infty)}$.
	
	\subsection{Core decomposition under $\widetilde{\mathbb{P}}^{(\infty)}$}
	\begin{proposition}[Core decomposition] \label{prop:coredecomposition} Under $ \widetilde{\mathbb{P}}^{(\infty)}$, the infinite core and the finite components dangling from it are independent and the latter are identically distributed with law $\mathbb{P}^{ {(\mathrm{free})}}$. 
	 \end{proposition}
	
	\begin{remark} In the finite gulp case, combining the previous proposition and the discussion at the end of the previous section, we deduce that under $ \mathbb{P}^{(\infty)}$, the only infinite simple component and the finite parts dangling form it are independent, and the latter are i.i.d.~of law $ \mathbb{P}^{( \mathrm{free})}$ except the component carrying the root edge which has law $ \mathbb{P}^{( \mathrm{free)}}$ biased by $1 + 2| \partial \map|$. This component might be reduced to $\dagger$ in which case the root edge belongs to the core.  See \cite{CMboundary} for a similar statement in the case of quadrangulations.
	\end{remark}
	
	 \noindent \textbf{Proof.} Inside $ \widetilde{\Map}^{(\infty)}$, denote by  $ C_{i}$, for $i \in \mathbb{Z}$, the finite components dangling from the infinite core, where $C_{0}$ is attached to the origin of the root edge. Fix $ j \geq 1$, and let us compute the law of $C_{-j}, C_{-j-1}, \dots,C_{0}$, $C_{1}, \dots, C_{j}$. 
	To do so, we explore $ \widetilde{\Map}^{(\infty)}$ with the metric exploration $ \mathcal{A}_{ \mathrm{metric}}$ until the minimal distance (inside the explored submap) from the origin to the exposed boundary is at least $j+1$. 
	
	 \begin{figure}[!h]
	  \begin{center}
	  \includegraphics[width=7cm]{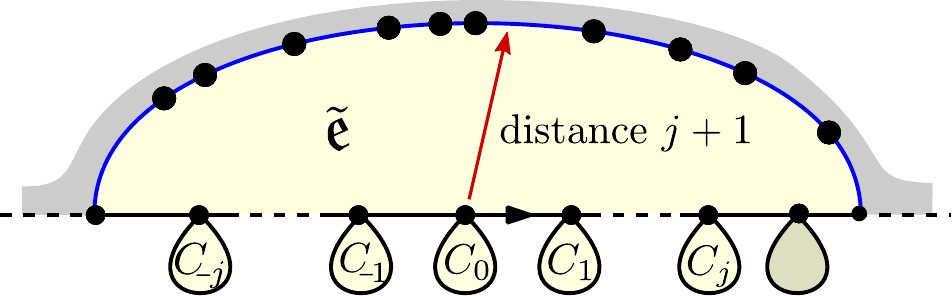}
	  \caption{ \label{fig:dangling} Schematic illustration of the result $\tilde{\mathfrak{e}}$ of a metric exploration. The parts $C_{-j},\ldots,C_j$ dangling from the core of $\tilde{\mathfrak{e}}$ agree with those dangling from the core of $\widetilde{\Map}^{(\infty)}$.}
	  \end{center}
	  \end{figure}
	  
	  The algorithm will almost surely stop and output an explored part $\tilde{\mathfrak{e}} \Subset \widetilde{\Map}^{(\infty)}$ which contains at least $j+1$ boundary edges of the core of $\widetilde{\Map}^{(\infty)}$ both on the left and on the right of the origin of the root edge. 
	  In particular, it has the same components $C_{-j},\ldots,C_j$ dangling from the core of $\tilde{\mathfrak{e}}$ (and possibly more further out to the left and right). The law of $\tilde{\mathfrak{e}}$ is characterized by  \eqref{eq:characSHP} and this can be turned into a product form where each of the component $C_{-j},\ldots,C_j$ contributes a weight
 $$ w_{ \mathbf{q}}(  C_i) \cdot  {c_{ \mathbf{q}}}^{ -| \partial C_i|}.$$ Recalling Section \ref{sec:freeboltz}, this means that $C_{-j},\ldots,C_j$ are i.i.d.~of law $ \mathbb{P}^{( \mathrm{free})}$.	  Since this is true for any $j$, the result follows. \qed \medskip 
	 
	\noindent \textbf{Proof of Proposition \ref{prop:spatialMarkovsimple}.} It suffices to prove the statement only for the region $ \hattilde{\mathfrak{e}}_{1}$ enclosed by the face incident to the left of the root edge, i.e.~the region discovered by a single simple-peeling step of the root edge. The full statement of the proposition follows by iteration of simple-peeling steps to discover $\hattilde{ \mathfrak{e}}$. Fix a map $ \mathbf{e}$. We use the definition of $\hat\Map^{(\infty)}$ as the core of $\widetilde{\Map}^{(\infty)}$ and scrutinize which events in $\widetilde{\Map}^{(\infty)}$ yield the  event  $\{\hattilde{\mathfrak{e}}_{1} =  \mathbf{e}\}$ in the proposition. If $2k$ is the degree of the face incident to the left of the root edge in $ \mathbf{e}$ then clearly the lazy-peeling of the root edge in $\widetilde{\Map}^{(\infty)}$ must discover a face $ \mathrm{f}$ of degree $2k$.
	According to Proposition \ref{prop:peelinghalfplaneMC} this happens with probability
	\begin{equation}\label{eq:fdegree}
	\widetilde{\mathbb{P}}^{(\infty)}( \mathrm{f} \text{ has degree }2k) = q_k c_{\mathbf{q}}^{k-1} H^\uparrow(2k-1).
	\end{equation}
	Furthermore, the remaining map $\widetilde{ \Map}^{(\infty)} \backslash  \mathrm{f}$ must have a core decomposition as in the following Figure \ref{fig:simple-peeling-onestep}. We will see that unless $k=1$ the map $\widetilde{ \Map}^{(\infty)} \backslash \mathrm{f}$ does \emph{not} have the same distribution as $\widetilde{\Map}^{(\infty)}$.
	 
	 \begin{figure}[!h]
	  \begin{center}
	  \includegraphics[width=1\linewidth]{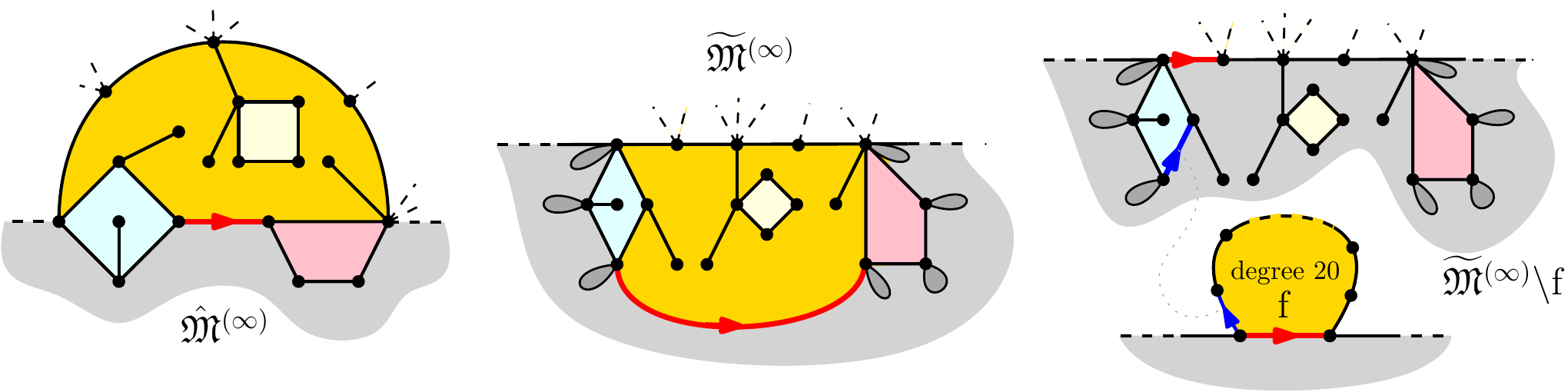}
	  \caption{ \label{fig:simple-peeling-onestep}On the left, the desired event for $\hat{\Map}^{(\infty)}$. In the center, the  corresponding possible events for $\widetilde{\Map}^{(\infty)}$ where the dark gray regions are arbitrary maps with a general boundary. On the right, we see the core decomposition of the map $\widetilde{\Map}^{(\infty)} \backslash  \mathrm{f}$. }
	  \end{center}
	  \end{figure}
	  
\begin{figure}[!h]
 \begin{center}
\includegraphics[width=13cm]{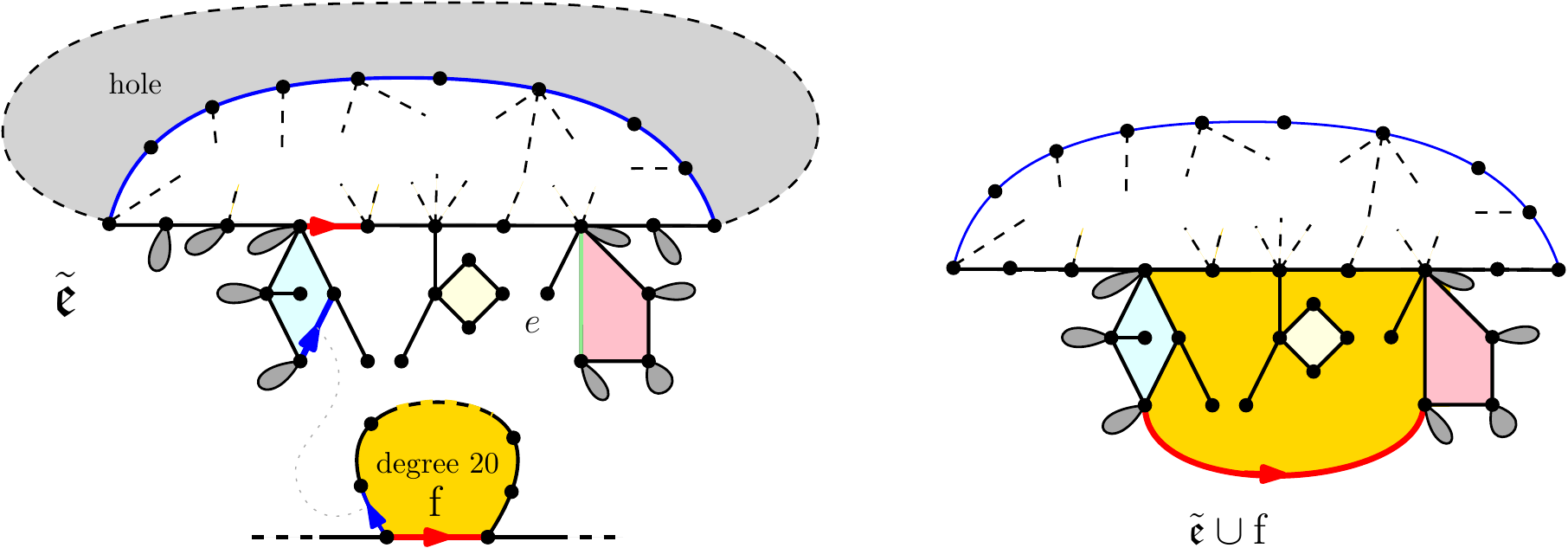}
 \caption{ \label{fig:setupef} Setup of the proof of the lemma. We fix a good submap $\tilde{ \mathfrak{e}}$ that can appear as a submap of $\Map[k]$, which means that after gluing a face of degree $2k$ on it (with proper starting and ending points) we get a possible submap of $\widetilde{\Map}^{(\infty)}$.}
 \end{center}
 \end{figure}

	  Let us write $\Map[k]$ as a short-hand for the map $\widetilde{ \Map}^{(\infty)} \backslash  \mathrm{f}$ conditionally on the degree of $\mathrm{f}$ being $2k$. If we distinguish the left-most oriented edge on the exposed boundary of $ \mathrm{f}$ this induces a blue edge on $\Map[k]$ (see Figure \ref{fig:simple-peeling-onestep}). We will root this map on  the first edge of the infinite core (in red in Figure \ref{fig:simple-peeling-onestep}) when tracing the contour of the face $ \mathrm{f}$ in $\widetilde{\Map}^{(\infty)}$. These two edges are equal or they have at most $2k-3$ edges in between them (since otherwise the root edge of $\widetilde{\Map}^{(\infty)}$ would not be on an infinite simple component). As announced above $\Map[k]$ is not distributed as $\widetilde{\Map}^{(\infty)}$ when $k \ne 1$ but the only difference lies in the map dangling from the origin of the red root edge:
	  
	\begin{lemma}\label{lem:Mkmarkov} The infinite core of $\Map[k]$ and the parts dangling from it are independent. The latter are identically distributed with law $ \mathbb{P}^{( \mathrm{free})}$ except for the component attached to the origin of the red root edge. This component $C_{ \mathrm{root}}$ has law 
	\begin{equation}\label{eq:Croot}
	\mathbb{P}( C_{ \mathrm{root}} = \map) = \frac{1}{\mathrm{W}_{c}\, H^{\uparrow}(2k-1)} w_{ \mathbf{q}}(\map) \cdot c_{ \mathbf{q}}^{-|\partial \map|} \left((2k-1) \wedge (2 | \partial \map|+1)\right),
	\end{equation}
	and conditionally on $C_{ \mathrm{root}}$ the blue edge of $\Map[k]$ is located uniformly either on the root edge or on one of the $(2k-2\wedge 2 |\partial C_{ \mathrm{root}}|)$ edges on the left of the root edge on $\partial C_{ \mathrm{root}}$. On the other hand, the core of $\Map[k]$ has law $\hat{ \mathbb{P}}^{(\infty)}$.
	\end{lemma}
	\noindent \textbf{Proof of the lemma.} We start by establishing a spatial Markov property for $\Map[k]$ with the help of Proposition \ref{prop:characSHP}. Let $\tilde{\mathfrak{e}} $ be a good submap with exposed boundary of length $p_{\mathrm{exposed}}$, equipped with another blue (oriented) edge on the component attached to the origin such that there are at most $2k-3$ between the blue edge and the root edge.
	Suppose furthermore that the edge $e$ to the right of the blue edge, chosen such that there are precisely $2k-3$ edges in between, is not adjacent to the hole of $\tilde{\mathfrak{e}} $.
	Then $\tilde{\mathfrak{e}}  \cup \mathrm{f}$ is the good submap obtained by adding to $\tilde{\mathfrak{e}} $ a new edge from the origin of the blue edge to the endpoint of $e$, and taking this new edge to be the root, see  Figure \ref{fig:setupef}.
	
By construction the exposed boundary of $\tilde{\mathfrak{e}}  \cup \mathrm{f}$ has length $p_{\mathrm{exposed}}$ as well.
	Then \eqref{eq:characSHP} and \eqref{eq:fdegree} imply that
	  \begin{eqnarray} \label{eq:loiMk} \mathbb{P}( \tilde{\mathfrak{e}}  \Subset \Map[k]) = \frac{ \widetilde{\mathbb{P}}^{(\infty)}( \tilde{\mathfrak{e}} \cup \mathrm{f} \Subset \widetilde{\Map}^{(\infty)})}{\widetilde{\mathbb{P}}^{(\infty)}( \mathrm{f} \text{ has degree }2k)} = w_{ \mathbf{q}}(  \tilde{\mathfrak{e}}) \cdot  {c_{ \mathbf{q}}}^{ | \partial^{*} \tilde{ \mathfrak{e}}| - | \partial \tilde{ \mathfrak{e}}|} \cdot \frac{H^\uparrow( p_{\mathrm{exposed}})}{H^\uparrow(2k-1)}.  \end{eqnarray}
	As in Proposition \ref{prop:coredecomposition}, we can write the above weight in a product form by isolating the contribution of the dangling parts and of the simple core of $ \tilde{ \mathfrak{e}}$. We deduce that the dangling parts and the core are independent, and that the core of $ \Map[k]$ has the same law as that of $\widetilde{\Map}^{(\infty)}$. The only difference with Proposition \ref{prop:coredecomposition} is that the component $C_{\mathrm{root}}$ at the root is size-biased by the number $(2k-1) \wedge (2|\partial C_{\mathrm{root}}|+1)$ of possible locations of the blue edge.
	The normalization constant in the denominator in \eqref{eq:Croot} is computed using
	$$\mathbb{E}^{( \mathrm{free})}[(2 |\partial \map|+1) \wedge (2k-1)] = \sum_{\ell=0}^{2k-2} \mathbb{P}^{{( \mathrm{free})}}(2|\partial \map| \geq \ell)\stackrel{\eqref{eq:Hdowndef}}{=} \sum_{\ell=0}^{2k-2} H^\downarrow(\ell) \stackrel{\eqref{eq:Hupdef}}{=} H^{\uparrow}(2k-1),$$
	which implies that
	$$ \sum_{\map} w_{ \mathbf{q}}(\map) \cdot c_{ \mathbf{q}}^{-|\partial \map|} \left((2k-1) \wedge (2 | \partial \map|+1)\right) \stackrel{\eqref{eq:Pfreedef}}{=} \mathrm{W}_c\,\mathbb{E}^{( \mathrm{free})}[(2 |\partial \map|+1) \wedge (2k-1)] = \mathrm{W}_c H^{\uparrow}(2k-1).$$
	\qed\medskip
	
	Coming back to the proof of the proposition, using the above core decomposition of $\Map[k]$ one can explicitly write down the probability inside $\widetilde{\Map}^{(\infty)}$ to produce an event yielding $\{\hattilde{\mathfrak{e}}_{1} =  \mathbf{e}\}$.
	Suppose the exposed boundary of $\mathbf{e}$ is of length $p$ and the (possibly empty) dangling components of $\mathbf{e}\setminus\mathrm{f}$ are $\tilde{C}_0, \tilde{C}_1, \ldots, \tilde{C}_p$.
	Let $\ell$ be the length of the outer boundary, i.e. the number of edges on the root face that are not shared with the hole, such that $p - \ell = 2( |\partial^* \mathbf{e}|-|\partial \mathbf{e}|)$. 
	Then
	\begin{align*}
	\hat{\mathbb{P}}^{(\infty)}\left(  \hattilde{\mathfrak{e}}_{1} = \mathbf{e} \right) &= \mathrm{W}_c^{\ell+1} \,\widetilde{\mathbb{P}}^{(\infty)}( \mathrm{f}\text{ has degree }2k)\, \widetilde{\mathbb{P}}^{(\infty)}( C_{\text{root}} = \tilde{C}_0) \prod_{i=1}^p \widetilde{\mathbb{P}}^{(\infty)}( C_i = \tilde{C}_i),	 
	\end{align*}
	where the factor of $\mathrm{W}_c^{\ell+1}$ takes into account the dangling maps in $\widetilde{\Map}^{(\infty)}$.
	With the help of \eqref{eq:fdegree} and Proposition \ref{prop:coredecomposition} this evaluates to
	\begin{align*}
	 \hat{\mathbb{P}}^{(\infty)}\left(  \hattilde{\mathfrak{e}}_{1}= \mathbf{e} \right)&= \mathrm{W}_c^{\ell+1}H^\uparrow(2k-1) \nu(k-1) \frac{1}{H^\uparrow(2k-1)} \prod_{i=0}^p \mathbb{P}^{(\text{free})}( \map = \tilde{C}_i )\\
	 &\stackrel{\eqref{eq:Pfreedef}}{=} \mathrm{W}_c^{\ell+1} w_{\mathbf{q}}(\mathbf{e}) c_{ \mathbf{q}}^{k-1} \mathrm{W}_c^{-p-1} \prod_{i=0}^p c_{ \mathbf{q}}^{-|\partial \tilde{C}_i|}.
	 \end{align*}
	 Finally, using that the length of the outer boundary of $\mathbf{e}\setminus \mathrm{f}$ is $p+\sum_{i=0}^p 2|\partial\tilde{C}_i|= \ell + 2k - 2$, we obtain
	 \begin{align*}
	 \hat{\mathbb{P}}^{(\infty)}\left(  \hattilde{\mathfrak{e}}_{1}= \mathbf{e} \right)= w_{\mathbf{q}}(\mathbf{e}) \left( \frac{c_{ \mathbf{q}}}{\mathrm{W}_c^2}\right)^{|\partial^*\mathbf{e}|-|\partial\mathbf{e}|} \stackrel{\eqref{def:cqhat}}{=} w_{\mathbf{q}}(\mathbf{e}) \hat{c}_{\mathbf{q}}^{|\partial^*\mathbf{e}|-|\partial\mathbf{e}|}.
	 \end{align*}
	 Conditionally on $\hattilde{\mathfrak{e}}_{1} = \mathbf{e}$ the remaining map is precisely the core of $\Map[k]$, which according to Lemma~\ref{lem:Mkmarkov} has the law of $\hat{\mathbb{P}}^{(\infty)}$.
	 This proves the proposition for the first step of a simple peeling exploration of $\hat{\Map}^{(\infty)}$.\qed  \medskip

	\subsection{Local limit of large boundaries}
	Now that the proper framework and properties are in place, deducing Theorem \ref{thm:localsimple} is straightforward. \medskip

	\noindent \textbf{Proof of Theorem \ref{thm:localsimple}.} Let $ \hattilde{ \mathfrak{e}}$ be a simple submap with one hole (of undetermined perimeter). Recall that $| \partial^* \hattilde{e}|$ and $| \partial \hattilde{e}|$ respectively denote the half perimeters of the exposed and internal boundaries. Then we have  	$$\hat{\mathbb{P}}^{(\ell)} (\hattilde{{\mathfrak{e}}} \hatSubset \map ) = w_{\mathbf{q}}(\hattilde{\mathfrak{e}}) \frac{\hat{W}^{(\ell+  | \partial \hattilde{e}|- | \partial^* \hattilde{e}|)}}{\hat{W}^{(\ell)}} \xrightarrow[\ell\to\infty]{\text{Prop. }\ref{prop:ratiosimple}} w_{\mathbf{q}}(\hattilde{\mathfrak{e}}) \hat{c}_{\mathbf{q}}^{ | \partial \hattilde{e}|- | \partial^* \hattilde{e}|} \stackrel{\text{Prop. }\ref{prop:spatialMarkovsimple}}{=} \hat{\mathbb{P}}^{(\infty)} ({\hattilde{\mathfrak{e}}} \hatSubset \map ).$$
	It is easy to see that this implies the weak convergence in the local topology.\qed  \medskip

 \section{Application to percolation}
 \label{sec:percosite}

Having defined our Boltzmann maps of the half-plane with a simple boundary, we can study Bernoulli percolation on them. We start by recasting the results of \cite{ACpercopeel} about bond and face percolation, in the context of triangulations and quadrangulations, to our more general setting of bipartite Boltzmann maps. We will then \emph{use} those two percolation processes to establish surprising connections between the lazy-peeling under $ \mathbb{P}^{(\infty)}$ and the simple-peeling under $ \hat{\mathbb{P}}^{(\infty)}$. Using those identities, we will be able to go through the argument of Richier \cite{Richier17b} and generalize its computation of the site-percolation thresholds to these random lattices.

\begin{center} \hrulefill \textit{ In this section we suppose that $ \gulp < \infty$.} \hrulefill  \end{center}

Thanks to the above hypothesis (which in particular implies criticality for $ \mathbf{q}$) we know (see the end of Section \ref{sec:Mtilde}) that the map $\hat{ \Map}^{(\infty)}$ of law $ \hat{ \mathbb{P}}^{(\infty)}$ can be defined as the unique infinite simple component of the map ${ \Map}^{(\infty)}$ of law $ { \mathbb{P}}^{(\infty)}$. In particular, since $ \Map^{(\infty)}\backslash \hat{\Map}^{(\infty)}$ is made of countably many finite components each separated from the infinite core by a vertex, the percolation thresholds (and the existence of an infinite cluster) are the same for percolations under $ \hat{ \mathbb{P}}^{(\infty)}$ and ${ \mathbb{P}}^{(\infty)}$. 

\begin{question} We leave as an open question to decide whether there is a non trivial phase transition for the (bond, face, site) Bernoulli percolation on $ \hat{ \Map}^{(\infty)}$ in the case $ \mathsf{g}_{ \mathbf{q}}= \infty$.
\end{question}

\subsection{Gulp and exposure}
\begin{figure}[!h]
 \begin{center}
 \includegraphics[width=9cm]{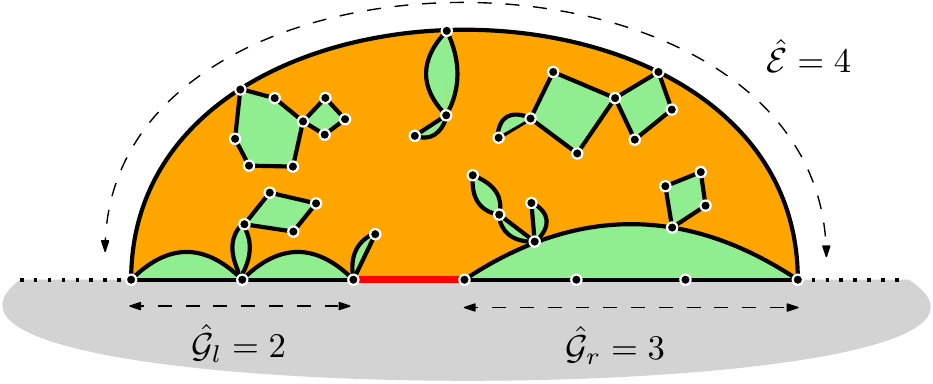}
 \caption{\label{fig:simplegulp}Example of a simple peeling step with $ \hat{\mathcal{E}}=4, \hat{\mathcal{G}}_{r}=3$ and $ \hat{\mathcal{G}}_{l}=2$. The light green regions  contain arbitrary maps with a simple boundary of the appropriate perimeter.}
 \end{center}
 \end{figure}

\begin{definition} Let the simple exposure $\hat{\mathcal{E}}$ and the simple gulps $\hat{\mathcal{G}}_{r}, \hat{\mathcal{G}}_{l}$ be the length of the outer perimeter and of the part of the inner perimeter on the right and left of the root edge in a one-step simple peeling transition under $\hat{\mathbb{P}}^{(\infty)}$. See Figure~\ref{fig:simplegulp}. The mean gulp and exposure are denoted by $ \sgulp =  \mathbb{E}[ \hat{\mathcal{G}}_{r}] = \mathbb{E}[ \hat{\mathcal{G}}_{l}]$ and $\sexpo=\mathbb{E}[ \hat{\mathcal{E}}]$.
\end{definition}

We recall also the analog notions in the case of lazy-peeling, where the gulp and exposure are easily expressed using the measure $\nu$ in \eqref{eq:nu}, see \cite[Definition 11.1]{CurStFlour}:
$$\expo = \sum_{k \geq 0} \nu(k) (2k+1), \quad \mbox{and } \quad \gulp = \sum_{k \geq 1}\nu(-k) (2k-1).$$
In the case of the lazy-peeling, \cite[Proposition 4.3]{BCMgasket} implies $\mathbf{q}$ is critical as soon as $\mathsf{g}_{ \mathbf{q}}< \infty$. By \cite[Proposition~4]{Bud15} the $\nu$-random walk must oscillate and we thus have $\expo = 2 \gulp+1$ or, in words,  the mean net change in the boundary length during one step of filled-in lazy-peeling under $ \mathbb{P}^{(\infty)}$ is zero. In light of this observation, it is natural to expect $\sexpo=2 \sgulp+1$, but we have no ``direct'' proof of this.
In fact, we will show the more surprising equality:

\begin{proposition} \label{prop:simplegulp} Suppose $\gulp < \infty$. Then the mean gulp and exposure coincide for the lazy and simple peeling process, i.e. $\gulp = \sgulp$ and $\expo= \sexpo$. In particular $ \sgulp, \sexpo$ are finite and $ \sexpo = 2 \sgulp+1$.
\end{proposition}

It is clear that in the dilute phase ($\gulp<\infty$), the mean simple exposure $\sexpo$ is finite, since conditionally on the root edge being on the core, the exposure in the first simple peeling step is bounded above by  the exposure in the lazy peeling step and $\expo < \infty$. But it is not clear \emph{a priori} that the mean simple gulp is finite! The proof of the above proposition is rather indirect since it uses the determination of the bond and face percolation thresholds using both the simple and the lazy peeling:

\noindent \textbf{Proof of Proposition \ref{prop:simplegulp}}. The almost sure bond and face percolation thresholds under $ \hat{ \mathbb{P}}^{(\infty)}$ have been determined in \cite{ACpercopeel} and are respectively equal to 
$$ \hat{p}_{c, \mathrm{bond}} = \frac{ \sgulp}{ \sgulp+1} \quad \mbox{and } \quad \hat{p}_{c, \mathrm{face}} =  \frac{\sgulp+1}{ \sexpo}.$$ 
Actually, \cite{ACpercopeel} only deals with the half-planar triangulation and quadrangulation with a simple boundary but the argument adapts readily, see \cite[Section 11.4.1]{CurStFlour} for a sketch. Similarly, the almost sure bond and face percolation thresholds under $ { \mathbb{P}}^{(\infty)}$ have been determined in \cite[Section 11.2 and 11.3]{CurStFlour} (see also \cite{BCMcauchy}) using the lazy-peeling process and are equal to 
$${p}_{c, \mathrm{bond}} = \frac{ \gulp}{ \gulp+1} \quad \mbox{and } \quad {p}_{c, \mathrm{face}} = \frac{ \gulp+1}{ 2\gulp+1}.$$
By the discussion at the beginning of this section, we know that the almost sure percolation thresholds in $\Map^{(\infty)}$ and $\hat{\Map}^{(\infty)}$ coincide. We can thus equate the last two displays  and deduce the proposition. 
\qed

 \subsection{Site percolation}
We now present the calculation of the site percolation threshold using an idea of Richier \cite{Ric15}. After sampling $  \hat{\Map}^{(\infty)} \sim \hat{\mathbb{P}}^{(\infty)}$ we color the vertices of the map independently in black with probability $p$ and white otherwise.
 
  \begin{theorem}[bond-percolation threshold on the half-plane]\noindent Suppose $ \mathbf{q}$ is an admissible weight sequence with $\gulp < \infty$ and put $$ p_{c, \mathrm{site}}= 1 - \frac{(\sum_{k=1}^\infty \nu(-k) )^2}{ 2 \nu(-1)
\gulp }.$$ Almost surely, there is no infinite black cluster for the $p$-Bernoulli site-percolation if $p \leq p_{c, \mathrm{site}}$ whereas there is an infinite black cluster if $p > p_{c, \mathrm{site}}$.
 \end{theorem}
 
 The proof of this theorem relies on the simple filled-in exploration method of \cite{Ric15} which we briefly recall. We speak of free-black-free boundary condition for a half-planar map if the vertices of the boundary of the map are all i.i.d.~black with probability $p$ and white otherwise, except for a finite connected segment of the boundary whose vertices are all black.

\begin{figure}
	\begin{center}
		\includegraphics[width=.6\linewidth]{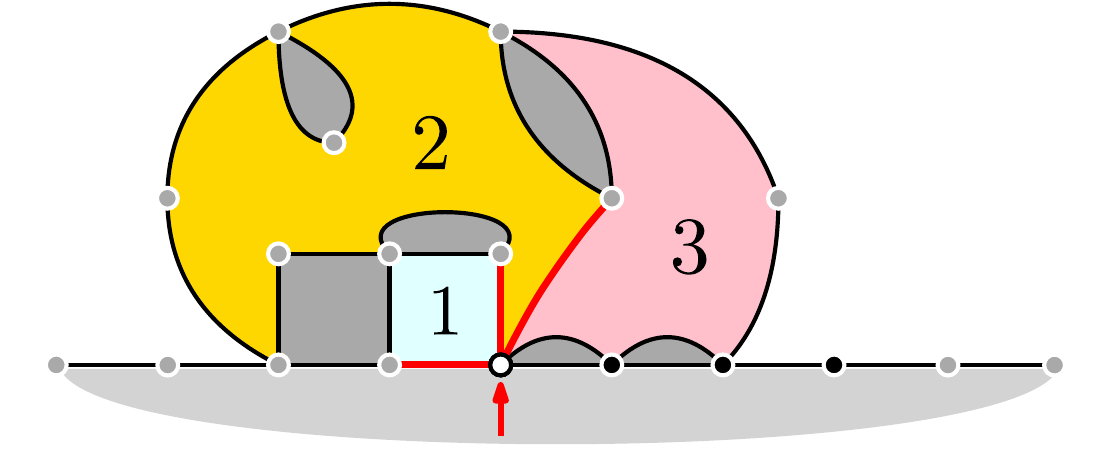}
	\end{center}
	\caption{\label{fig:peelingvertex}Simple peeling of a white vertex: when a white vertex is discovered, we peel it by revealing the faces (3 in our example) clockwise around it until the vertex is swallowed. The dark gray holes are filled-in with $\partial$-simple maps of the appropriate perimeter. The gray vertices are free, i.e.~independently sampled to be black with probability $p$ or white with probability $1-p$.}
\end{figure}

\begin{center}
\fbox{\begin{minipage}{15cm} \textbf{Algorithm $ \mathcal{A}_{ \mathrm{site}}$:} Suppose $ \hattilde{\mathfrak{e}}_{n}$ has a free-black-free boundary condition. The algorithm first reveals the color of the first free vertex to the left of the black boundary. If it is black then we move on to the next step. If it is white, then we peel this vertex. To do so, we iteratively peel (via simple peeling steps) the edge just to the left of this white vertex until we encounter a simple peeling step that swallows it, i.e.~such that $ \hat{ \mathcal{G}}_{r}>0$. We also fill in all the finite holes encountered along the way. See Figure \ref{fig:peelingvertex}.
\end{minipage}}
\end{center}
 
It is easy to check that the free-black-free boundary condition is preserved during such an exploration and that as long as the length of the black boundary (i.e.~number of black vertices) is positive it evolves as a random walk killed on $ \mathbb{Z}_{<0}$ with i.i.d.~increments distributed as 
 $$ \epsilon - (1- \epsilon) \cdot \left((\hat{ \mathcal{G}}_{r} \mid  \hat{ \mathcal{G}}_{r}>0) -1\right),$$ where $\epsilon$ is a Bernoulli random variable with success parameter $p$ independent of $(\hat{ \mathcal{G}}_{r} \mid  \hat{ \mathcal{G}}_{r}>0)$, which has the law of the simple gulp conditioned to be strictly positive. One can easily adapt \cite[Lemma 3.5 and Theorem 1.1]{Ric15} to our context and get the statement of the above theorem where $p_{c, \mathrm{site}}$ is determined by the requirement that the above random walk has zero drift, i.e.
 \begin{eqnarray} \label{eq:pcsite} p_{c, \mathrm{site}} -(1 - p_{c, \mathrm{site}}) \mathbb{E}[ \hat{ \mathcal{G}}_{r}-1 \mid \hat{\mathcal{G}_{r}}>0] =0. \end{eqnarray}
 Hence we have
 $$ p_{c, \mathrm{site}} \underset{ \eqref{eq:pcsite}}{=}  1 - \frac{\mathbb{P}( \hat{ \mathcal{G}}_{r}>0)}{\mathbb{E}[ \hat{ \mathcal{G}}_{r}]} \underset{ \mathrm{Prop.} \ref{prop:simplegulp}}{=} 1-  \frac{\mathbb{P}( \hat{ \mathcal{G}}_{r}>0)}{\gulp},$$
and  our theorem is thus completed by the following calculation: 
\begin{lemma}\label{lem:gulpcalc} We have 
$$ \mathbb{P}( \hat{ \mathcal{G}}_{r}>0) = \frac{(\sum_{k=1}^\infty \nu(-k) )^2}{ 2 \nu(-1)}$$
\end{lemma}
Notice that in the notation of Section \ref{sec:simpleenu} we have $\mathbb{P}( \hat{ \mathcal{G}}_{r}>0) = 1 / \hat{c}_ \mathbf{q}$, which is begging for a simpler explanation than the proof we give below.\medskip 

\noindent \textbf{Proof.} 
Using the same strategy and notation as in the proof of Proposition \ref{prop:spatialMarkovsimple}, we examine the one-step simple peeling of $\hat{\Map}^{(\infty)}$ by looking at the one-step lazy peeling of $\widetilde{\Map}^{(\infty)}$.
Conditionally on revealing a face of degree $2k$, which happens with probability \eqref{eq:fdegree}, the event $\mathcal{G}_{l}=0$ corresponds precisely to the blue edge of the root component $C_{\mathrm{root}}$ having its origin incident to the infinite core of $\Map[k]$,
\begin{align*}
\mathbb{P}( \mathcal{G}_{l}=0) &= \sum_{k=1}^\infty \widetilde{\mathbb{P}}^{(\infty)}(\text{f has degree }2k)\, \mathbb{P}(\text{blue root of }C_{\mathrm{root}}\text{ in }\Map[k]\text{ is incident to infinite core})\\
&=\sum_{k=1}^\infty q_k c_{\mathbf{q}}^{k-1} H^\uparrow(2k-1)\,  \mathbb{P}(\text{blue root of }C_{\mathrm{root}}\text{ in }\Map[k]\text{ is incident to infinite core}).
\end{align*}
\begin{figure}[t]
	\begin{center}
		\includegraphics[width=.9\linewidth]{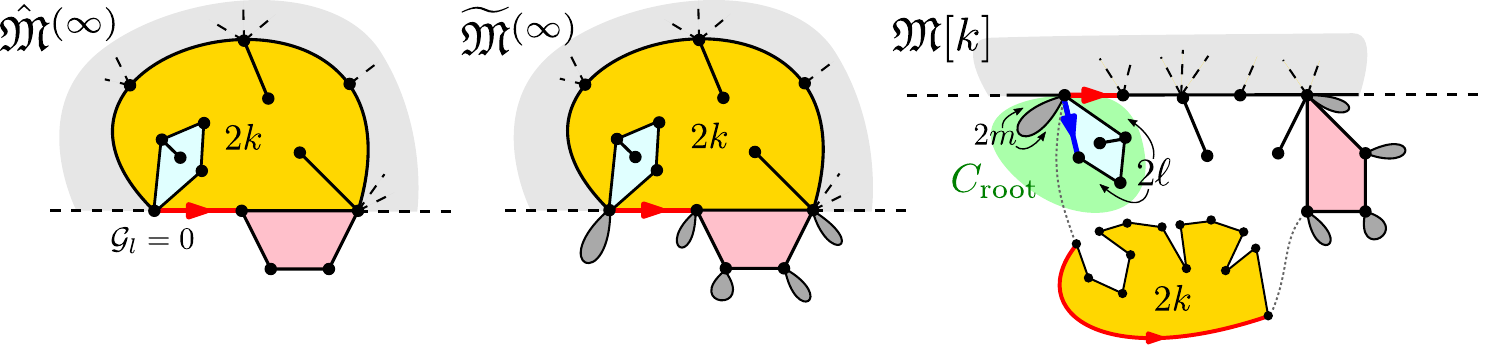}
	\end{center}
	\caption{\label{fig:gulpzero} We examine the situation in $ \widetilde{\Map}^{(\infty)}$ when the event $\mathcal{G}_l=0$ during the peeling of a face of degree $2k$ happens in $\hat{ \Map}^{(\infty)}$. If we remove the face of degree $2k$, the remaining map $\Map[k]$ must have the property that the dangling component $C_{\mathrm{root}}$ at the origin (shaded in green) decomposes into two free Boltzmann maps (of perimeters $2m$ and $2\ell$ respectively). }
\end{figure}%
Using the explicit law of $C_{\mathrm{root}}$ in Lemma \ref{lem:Mkmarkov} and its decomposition into two free Boltzmann maps of perimeters $2m$ and $2\ell$ (right picture in Figure \ref{fig:gulpzero}), we find
\begin{align*}
\mathbb{P}( \mathcal{G}_{l}=0) &= \sum_{k=1}^\infty \widetilde{\mathbb{P}}^{(\infty)}(\text{f has degree }2k)\, \mathbb{P}(\text{blue root of }C_{\mathrm{root}}\text{ in }\Map[k]\text{ is incident to infinite core})\\
& = \frac{1}{\mathrm{W}_c}\sum_{k=1}^\infty q_k c_{\mathbf{q}}^{k-1} \sum_{\ell=0}^{k-1} W^{(\ell)} c_{\mathbf{q}}^{-\ell} \sum_{m=0}^\infty W^{(m)}c_{\mathbf{q}}^{-m} 
=\sum_{k=1}^\infty q_k c_{\mathbf{q}}^{k-1} \sum_{\ell=0}^{k-1} W^{(\ell)} c_{\mathbf{q}}^{-\ell}.
\end{align*}
Expressed in terms of the measure $\nu$ of \eqref{eq:nu}, which in particular satisfies $\nu(-1) = 2 / c_{\mathbf{q}}$, this becomes
\[
\mathbb{P}( \mathcal{G}_{l}=0) = \frac{1}{\nu(-1)} \sum_{k=1}^\infty \nu(k-1) \sum_{\ell=0}^{k-1}\nu(-\ell-1).
\]
The sum on the right-hand side can be further simplified using that $\nu$ is normalized, 
\[
1 = \left(\sum_{k=-\infty}^\infty \nu(k)\right)^2 = \left(\sum_{k=0}^\infty \nu(k)\right)^2 + \sum_{\ell=1}^\infty\sum_{k=-\infty}^\infty \nu(k)\nu(-\ell-k-1) + 2 \sum_{k=1}^\infty \nu(k-1) \sum_{\ell=0}^{k-1}\nu(-\ell-1),
\]
and satisfies Tutte's equation (see e.g.\ \cite[Eq. (9)]{Bud15}) rewritten in terms of $\nu$ as 
\[
\nu(-\ell-1) = \frac{1}{2} \sum_{k=-\infty}^\infty \nu(k)\nu(-\ell-k-1).
\]
Together these identities imply that
\begin{align*}
\mathbb{P}( \mathcal{G}_{l}=0) & = \frac{1}{\nu(-1)} \left[ \frac{1}{2}-\frac{1}{2}\left(1 - \sum_{\ell=1}^{\infty} \nu(-\ell)\right)^2 - \sum_{\ell=1}^{\infty} \nu(-\ell-1)\right]= 1 - \frac{\left( \sum_{\ell=1}^{\infty} \nu(-\ell)\right)^2}{2\nu(-1)}
\end{align*}
and the result follows from the fact that $\mathcal{G}_{l}$ and $\mathcal{G}_{r}$ are equal in distribution and non-negative.
\qed  
	

\end{document}